\documentclass[1p,final]{article}
\usepackage{graphicx}
\usepackage[english]{babel}
\usepackage{amsmath,amssymb,amsfonts}
\usepackage{setspace}
\newcommand{\obj}[1]{\mathbf{#1}}
\newcommand{\objf}{\obj{f}}
\newcommand{\objh}{\obj{h}}
\newcommand{\objdist}[2]{\lVert #1-#2 \rVert}

\begin{document}

\title{A~New Mechanism for Maintaining Diversity of Pareto Archive in Multiobjective Optimization}

\author{Jaroslav H\' ajek$^1$,
Andr\' as Sz\" oll\" os$^1$, 
Jakub \v S\'\i stek$^2$\dag \\ 
{\small $^1$ Aeronautical Research and Test Institute,} \\
{\small Beranov\' ych 130, Praha, CZ-199 05, Czech Republic} \\
{\small $^2$ Institute of Mathematics of the AS CR,} \\
{\small \v Zitn\'a 25, Praha, CZ-115 67, Czech Republic} \\
{\small \dag Corresponding author} \\
{\small \tt hajek@vzlu.cz, szollos@vzlu.cz, sistek@math.cas.cz} 
}

\date{}

\maketitle

\begin{abstract}
The article introduces a~new mechanism for selecting individuals to a~Pareto archive.
It was combined with a~micro-genetic algorithm and tested on several problems. 
The ability of this approach to produce individuals uniformly distributed along the Pareto set
without negative impact on convergence is demonstrated on presented results.
The new concept was confronted with NSGA-II, SPEA2, and IBEA algorithms from the PISA package.
Another studied effect is the size of population versus number of generations for small populations.

\end{abstract}

\section{Introduction}
\label{sec:intro}

The field of multiobjective design optimization has evolved very fast during last years,
reflecting the need of solving tasks with several conflicting criteria, 
which is common in practical problems.
From the mathematical point of view,
this corresponds to minimization/maximization of a~vector-valued function,
which rarely leads to a single solution.
Consequently,
a whole hyperplane of trade-off solutions, called Pareto-optimal set, is expected as the result instead of a~single optimum.

A number of algorithms have been presented that generate a~set of solutions approximating this hyperplane.
The quality of the approximation is usually considered from two points of view:
(i) the closeness to the exact trade-off surface and
(ii) its distribution.
The former is related to convergence properties of an algorithm
while the latter describes its ability to maintain diversity.
An ideal algorithm should produce well converged solutions perfectly distributed along the Pareto front.
However, these requirements are conflicting, and many current approaches concentrates on one of them 
finding reasonable compromise in the other.

In this study, 
our attention is focused on the second aspect of diversity of the Pareto-optimal set,
namely we present a~new strategy for maintaining variety of members of a~Pareto archive.

The problem of maintaining uniform distribution at an affordable cost has been addressed by many algorithms.
It is known that the notion of crowding distance proposed by Deb et al. for algorithm NSGA-II \cite{Deb-2001-MOO,Deb-2000-FEN} 
is not sufficient to maintain diversity of the evolution for more than two objectives (e.g. \cite{Deb-2003-FME,Deb-2002-SMO}).
On the other hand, SPEA2 by Zitzler et al. \cite{Zitzler-2002-EMD} is usually able to produce well spread solutions even for three or more objectives.
The concept of archiving promising design vectors was first introduced for SPEA by Zitzler and Thiele \cite{Zitzler-1999-MEA}.
Knowles and Corne presented the Pareto Archived Evolutionary Strategy (PAES) \cite{Knowles-1999-PAE} and proposed the adaptive grid algorithm (\cite{Knowles-2000-ANF}) to maintain diversity.
However, it is difficult to keep the efficiency of this approach in cases with more than three objectives.

The new mechanism presented in this paper was implemented in micro-genetic algorithm $\mu$ARMOGA proposed by Sz\"oll\"os et al. \cite{Szollos-2008-AOM},
and results for three standard three-objective benchmark problems are presented.

Our second aim is to investigate the effect of population size for small (sometimes called micro) populations on the performance of $\mu$ARMOGA.
It was reported by Krishnakumar \cite{Krishnakumar-1989-MGA} for single-objective optimization and by Coello and Pulido \cite{Coello-2001-MGA} for multiple objectives, 
that very small populations can lead to fast convergence to the Pareto front.
In this context, most experiments were performed using populations of 4, 10 and 20 individuals.

Results got by $\mu$ARMOGA equipped with the new archiving mechanism are compared with those obtained by two leading methods in the field,
namely NSGA-II by Deb et al. \cite{Deb-2000-FEN} and SPEA2 by Zitzler et al. \cite{Zitzler-2002-EMD}, 
and a~recent interesting algorithm IBEA by Zitzler and K\"unzli \cite{Zitzler-2004-IBS}.
All these are implemented in the Platform and Programming Language Independent Interface for Search Algorithms (PISA) 
\footnote{\texttt{http://www.tik.ee.ethz.ch/sop/pisa}} \cite{Bleuler-2003-PIS}.
PISA is an~interesting open source package developed by the team of Prof. E. Zitzler at ETH Z\" urich.
The software implements various selection, crossover, and mutation operators and objective function evaluations.
An important idea of the project is to separate the selection of promising candidates from objective function evaluation, 
crossover and mutation and implement these in two separate programs, interchanging information via formatted files.
These programs are called \emph{selectors} and \emph{variators} in the PISA context.
There is an~increasing number of ready-to-use variators and selectors that can be downloaded from the web page of the PISA project. 
Therefore, the system offers a~simple way to produce fair comparisons of various selection schemes with the same variator.
While the described scheme of splitting an~evolutionary algorithm into two separate programs is very useful for some techniques, 
in our opinion, it does not fit to algorithms with strong coupling between both stages via the use of an archiving procedure. 
That is the reason why our implementation of $\mu$ARMOGA was used, 
instead of integrating the proposed archiving technique into the PISA framework.

Three metrics, 
measuring both convergence to the exact front and diversity of the approximate set,
are used for the comparison.
It is observed, that the new algorithm produces very good distribution of individuals
outperforming in this respect the other algorithms in many cases.
The archiving strategy does not seem to affect its convergence.
Moreover, diversity is maintained in an affordable way as suggested by presented numerical experiments.

The rest of the paper is organised as follows.
In Section~\ref{sec:armoga}, $\mu$ARMOGA is recalled with an~emphasis on its main aspects.
Section~\ref{sec:archive} contains the main contribution,
which is the proposition of a new archiving mechanism.
Tests and comparisons with the other evolutionary techniques can be found in Section~\ref{sec:comparison},
where we describe the test problems (\ref{sec:problems}),
metrics used for evaluating the performance (\ref{sec:metrics}),
detailed setting of particular algorithms (\ref{sec:setting}),
and organization of the experiments (\ref{sec:experiments}), respectively.
Our findings are discussed in detail in Section~\ref{sec:discussion}, 
while Section~\ref{sec:conclusion} contains summary of the work and concluding remarks.

\section{$\mu$ARMOGA -- Multiobjective micro-genetic algorithm with range adaptation: A Brief Introduction}
\label{sec:armoga}

To minimize the costly evaluation of individuals, 
it is straightforward to see that one way to go is to minimize their number. 
It is well known for evolutionary practitioners,
that using smaller populations and applying the evolutionary operators many times is often more favourable 
than vice versa (e.g. \cite{Oyama-2000-WDU}). 
This idea can be brought to an~extreme by using a micro-population (e.g. 4, 5, 10 individuals), 
what we really did when we utilized some ideas of Krishnakumar \cite{Krishnakumar-1989-MGA} and of Coello and Pulido \cite{Coello-2001-MGA}.
Krishnakumar came with the concept of micro-genetic approach first, and used it for single-objective optimization. 
His algorithm contained only selection and crossover operators, and no mutation operator. 
Instead, the author introduced a reinitialization technique, which was invoked once in a~few generations to ensure diversity for the evolution. 
The latter two researchers proposed a~micro-genetic algorithm enabling to tackle multi-objective problems. 
Their concept was similar to that of Krishnakumar, i.e. it contained selection, crossover, and reinitialization operators 
supplemented by a mutation operator. 
Both algorithms were verified on various test problems. 
In both the cases, the micro-genetic variants converged to the optimum (Pareto-front) much faster 
than their macrogenetic counterparts used for comparison.

In the approach by Sz\"oll\"os \cite{Szollos-2008-AOM}, 
microgenetic algorithm is supplemented by range adaptation and ``knowledge-based'' reinitializion procedure exploiting the Pareto-archive to generate better individuals.

The concept of range adaptation was originally introduced by Arakawa and Hagiwara~\cite{Arakawa-1998-DAR}, 
who used it with binary coding of the design variables. 
Its essence lies in ability to promote the evolution towards promising regions of the design space via sophisticated manipulation with the population statistics.
Due to the coding, it contained some artificial parameters which were hard to guess in general. 
Oyama \cite{Oyama-2000-WDU} used range adaptation in real domain and was successful in avoiding this drawback via encoding a design variable to a real number $r_i \in (0,1)$ 
defined by integration of the Gaussian distribution $N(0,1)$ 
\begin{equation}
r_i = \int_{-\infty}^{\widetilde{p}_i}N(0,1)(z)\mbox{d}z,
\end{equation}
where $\widetilde{p}_i$ is linked to the original design variable $p_i$ by
\begin{equation}
p_i = \sigma _i \cdot \widetilde{p}_i + \mu _i.
\end{equation}
We are using this encoding scheme too, with one important difference: 
Oyama originally calculated the average $\mu _i$ and the standard deviation $\sigma _i$ by sampling the upper half of the population, 
which is justified as long as macro-populations are used (e.g. with more than fifty individuals). 
But such approach would be too restrictive in the case of microevolution, 
since the upper half of the micropopulation contains too little information to keep the diversity. 
Consequently, the evolution quickly ends up in premature convergence. 
Thus, we calculate both by taking into account the whole population.

\begin{figure}[p]
\centering
\includegraphics[height=70mm]{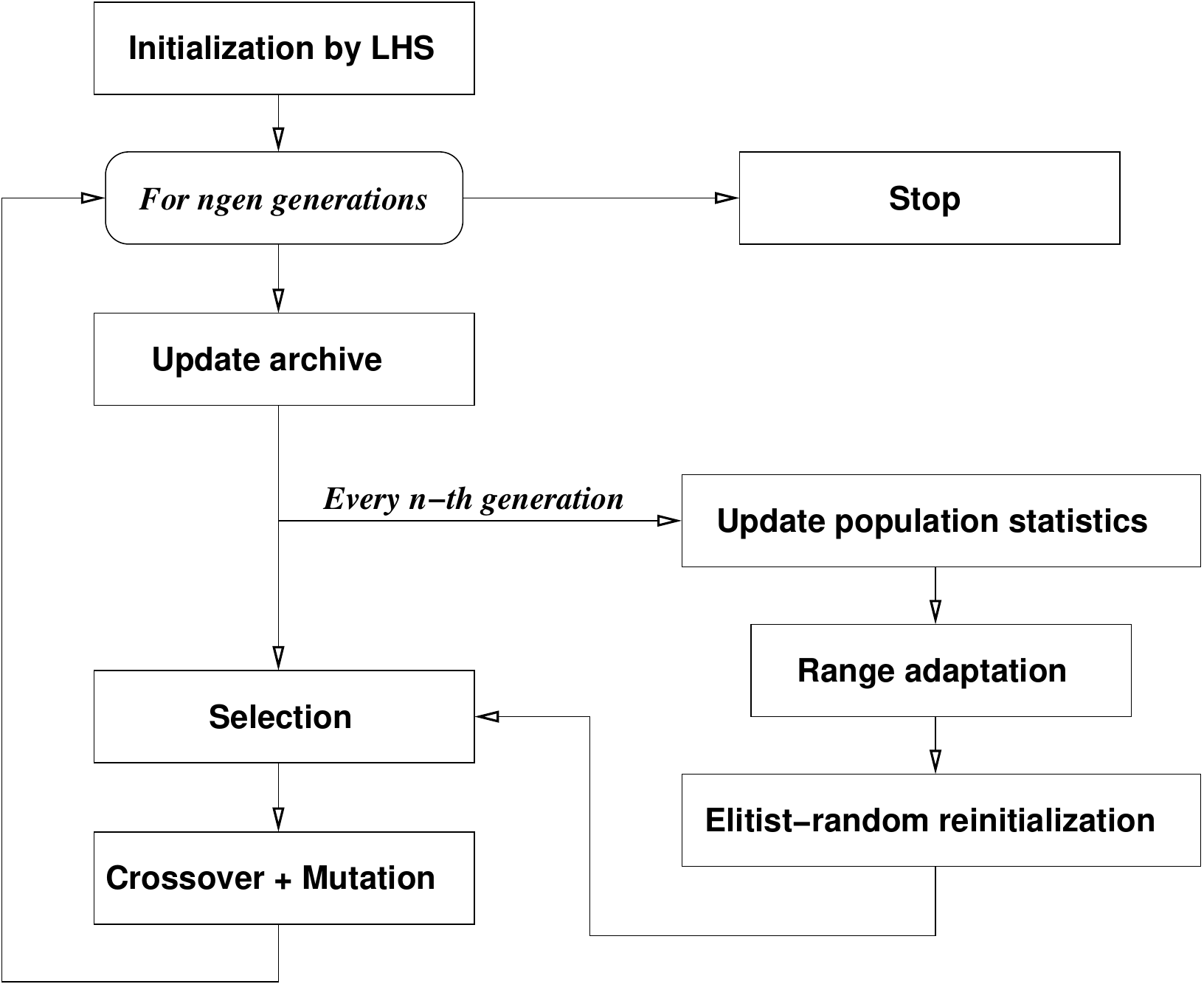}
\caption{\label{fig:miarmoga_scheme} Simple scheme of the $\mu$ARMOGA algorithm.}
\end{figure}

``Knowledge based'' reinitialization resulted from an attempt to use the members of the Pareto-archive to get new members, 
superseding the old ones by putting several of them into the reinitialized population. 
Moreover, only a subset of the archive is considered. 
For instance, two archive members with extreme values of two different objectives chosen randomly are usually exploited. 
In this way, it is possible to further improve the whole archive by improving its subsets. 
The functioning of $\mu$ARMOGA can be seen in Figure \ref{fig:miarmoga_scheme}.
After initialization of the population by Latin hypercube sampling (LHS) and evaluation depicted as archive update, 
the evolution goes through selection, mating and mutation to evaluation of the new population. 
Each $n$-th generation the population statistics is updated, 
range-adaptation takes place, followed by knowledge based (elitist-random) reinitialization. 
A~thorough description of the algorithm is to be found in \cite{Szollos-2008-AOM}. 

Our approach contains two new system parameters: adaptation factor $\delta$ and minimal standard deviation $\sigma_{min}$. 
In short, $\mu$ARMOGA strives to keep the evolution in a~permanently ``excited'' state via forced modification of the population statistics. 
It practically means that the standard deviation is not allowed to fall under certain minimal value of $\sigma_{min}$ for any design variable. 
This helps to prevent the micro-genetic algorithm from getting stuck in premature convergence. 
The role of the adaptation factor $\delta$ lies in controlling the frequency of range adaptation: if reinitialization is necessary, 
and the new standard deviation of a~design variable is changed by more than $\delta \cdot \sigma_{old}$, 
where $\sigma_{old}$ is the standard deviation when the last reinitialization took place, 
then the range of that design variable is adapted.

\section{The archiving algorithm}
\label{sec:archive}

\textit{Pareto archive} is a~key component of many evolutionary algorithms.
It acts as a~collector of good individuals during the evolution, and is often
used to give the resulting Pareto front approximation at the end of the evolution.
After new individuals are evaluated, the archive is improved if these individuals dominate or are non-dominated with respect to 
the existing individuals of the archive. 
During reinitialization,
the micro-genetic algorithm retrieves information
from the archive, using it to explore the promising regions of the search
space. 

Obviously, in any real setup we must limit the number of individuals stored in an archive.
This is necessary not just to keep the amount of information processed feasible, 
but also to get a~good diversity of the resulting approximation.
In our strategy, we use a~fixed upper limit on the number of individuals stored
in the archive. Ideally, we want to end up with a~full archive of Pareto-optimal
solutions that is ``well spread'' over the true Pareto front of the problem.
Our approach is an~archive dealing with a~single new individual at a~time. 
This is particularly suitable for micro-evolutionary approaches, where we only
have a few new individuals from each generation.

When a~new individual arrives, it is first checked for Pareto dominance with all existing
members of the archive. Now, we distinguish among three cases:
\begin{itemize}
\item The new individual is dominated by one or more members of the archive.
In this case, the new individual is discarded.
\item The new individual dominates one or more members of the archive.
The dominated ones are removed, and the new individual is added to the
archive and the internal information of the archive is updated (see below).
\item The new individual is non-dominated and non-dominating. If the number of
members of the archive has not yet reached the upper limit, the new
individual is added as in the previous case. In the opposite case, we need to
discard at least one individual (either the newcomer or one from the archive),
but we can not decide this by Pareto dominance. 
In this case, we proceed to 
the secondary decision procedure described below:
\end{itemize}

If we arrive at the case that can not be resolved by Pareto dominance, our secondary
goal is to maximize the distance between neighbouring individuals, 
based on some distance-measure in the objective space. 
In this paper, we use the standard Euclidean distance, 
which is meaningful for any dimension of the objective space.

First, we consider the minimum pairwise distance, i.e.,
\begin{equation}\label{eq:pwmin}
\min_{\substack{i,j \in P \\ i\neq j}} \objdist{\objf_i}{\objf_j},
\end{equation}
where $P$ denotes the set of archived individuals and $\objf_i$ stands
for the vector of objective values of individual $i$.
We take the pair of individuals that achieves the minimum in the above expression.
If there are multiple pairs, we take any of them. Without the loss of generality, we
assume that the minimum pair is $\objf_1, \objf_2$. Further, we denote 
the vector of objective values of the new individual as $\objh$. 
If 
\begin{equation}\label{eq:globimp1}
\min_{\substack{k \in P \\ k\neq 1}}\objdist{\objf_k}{\objh} > \objdist{\objf_1}{\objf_2},
\end{equation}
we can replace $\objf_1$ by $\objh$.
Alternatively, if
\begin{equation}\label{eq:globimp2}
\min_{\substack{k \in P \\ k\neq 2}}\objdist{\objf_k}{\objh} > \objdist{\objf_1}{\objf_2},
\end{equation}
we can replace $\objf_2$ by $\objh$.
If either of the above conditions is satisfied, the overall minimum pairwise distance
will be improved by the substitution or, if there were multiple minimal pairs, it will
stay the same but the number of minimal pairs will reduce. We call this as the 
\textit{global improvement} check.

If neither of these conditions is satisfied, we consider the closest
archived individual to $\objh$, say, $\objf_c$ instead. If 
\begin{equation}\label{eq:locimp}
\min_{\substack{k \in P \\ k\neq c}}\objdist{\objf_k}{\objh} > \min_{\substack{k \in P \\ k\neq c}}\objdist{\objf_k}{\objf_c},
\end{equation}
we replace $\objf_c$ by $\objh$. If this condition holds, there is a~certain
subset of the archived individuals whose pairwise minimum will improve.
This is the \textit{local improvement} check.

If neither check is successful, we discard the new individual.

Searching for the minimum-distance pair of the archive afresh each time an
individual is considered would be too costly. To make the procedure efficient,
we maintain for each archived individual a~pointer to its closest neighbour
(or any of them). Therefore, searching for the pairwise minimum in \eqref{eq:pwmin}
requires only one pass through the archive. Similarly, the right-hand side of
equation \eqref{eq:locimp} is simply the distance of $\objf_c$ to its closest
neighbour. Hence, these two checks only require computing the distances of the
new individual to all archived individuals, and computing the minima on left-hand
sides of the equations \eqref{eq:globimp1}, \eqref{eq:globimp2}, and \eqref{eq:locimp}.
Thus, \textit{deciding} whether to add a~new individual has linear complexity in terms
of number of archived individuals (evaluating mutual pairwise dominance also has
linear complexity).

If the new individual is to be added, the existing closest-neighbour links need
to be updated. Each resulting archive member is considered in turn.  If the
link is valid (i.e. the closest neighbour in the archive was not discarded),
we simply check if the newcomer is closer, and possibly update the link.
This takes only constant time. However, if the link became invalid (the former
closest neighbour was discarded), we need to compute the closest neighbour
afresh by computing objective distances of the updated individual to all
others. 

It can be proven by a simple argument based on $k$-dimensional ball volumes
that the maximum number of points in $k$-dimensional space having a single
common closest neighbour is bounded from above by a constant depending on $k$.
Since the Pareto archive consists of mutually non-dominating vectors, which
can not be arranged arbitrarily, in our case the constant is even smaller.
For instance, for a two-objective optimization, i.e. $k=2$, a single archive
member can be the closest neighbour to at most two other members at the same
time. 

Using this argument, it can be easily seen that the complexity of a single
archive update has complexity $O(N + MN)$, where $N$ is the size of the Pareto
archive and $M$ is the number of archive members dominated by the new individual.
As was already said, merely \textit{deciding} whether the newcomer is to be
added costs $O(N)$. If the decision is positive, there are two cases:
either the newcomer dominates some $M$ existing archive members, or it was added
based on the secondary decision procedure.
In the former case, $M$ members will be discarded, so at most $c\,M$ nearest-neighbour
links will need to be updated, $c$ being the upper bound constant discussed in the
previous paragraph. In the latter case,
one existing member is discarded, so at most $c$ existing links must be updated.
Given that updating a single link costs $O(N)$, together we have the cost
$O(N+c(1+M)N)$ which can be simplified to $O(N+MN)$, given that $c$ is a constant
independent of $M,N$. 

While in principle $M$ can be as high as $N$, in practice it drops to $M \ll N$ very
quickly as the convergence proceeds and new dominating individuals become increasingly
rare. It should also be noted that if $M$ is high at one step, the evolution
continues with an archive of $N-M$ which will be significantly smaller than $N$.
Numerical experiments confirm that in real evolutionary runs, the average number
of invalid links per archive update is very small, even much smaller than the 
theoretical bounds suggested above. 
This might be observed from Tables~\ref{tab:updates} and \ref{tab:asize_updates}. 
Therefore, we can conclude that the procedure of adding new individual to our
archive is essentially of linear complexity.

\section{Comparison of results}
\label{sec:comparison}

The abilities of the new archiving mechanism are first demonstrated on test functions DTLZ1, DTLZ2 and DTLZ4,
suggested by Deb et al. \cite{Deb-2002-SMO}.
To examine the influence of population size, our algorithm was run separately with 4, 10 and 20 individuals.
Obviously, it is preferable to maintain the number of function evaluations as low as possible, therefore
we study the behaviour of the aforementioned approaches for three fixed numbers of evaluations, 
4\,000, 20\,000 and 40\,000.
For test problem DTLZ1, the number of function evaluations is extended to 100\,000 and 200\,000,
since the algorithms were unable to converge to the global Pareto front with just 40\,000 computations.

To further investigate the behaviour of the proposed method, we performed an experiment with test problem WFG1 
suggested by Huband et al. \cite{Huband-2005-SMO}.
For this difficult problem, it was necessary to run the evolution to as many as 2\,000\,000 evaluations to obtain
reasonable convergence to the Pareto front.

\subsection{Test problems}
\label{sec:problems}

The algorithms are compared on three benchmark problems introduced in \cite{Deb-2002-SMO}.
The following form of them is considered:
\begin{itemize}
\item DTLZ1

Minimize $f_1$,$f_2$,$f_3$, where
\begin{eqnarray}
f_1(\mathbf{x}) &=& \frac{1}{2} x_1 x_2 (1 + g(\mathbf{x}_M)), \\
f_2(\mathbf{x}) &=& \frac{1}{2} x_1 (1 - x_2) (1 + g(\mathbf{x}_M)), \\
f_3(\mathbf{x}) &=& \frac{1}{2} (1 - x_1) (1 + g(\mathbf{x}_M)), \\
g(\mathbf{x}_M) &=& 100 \left( 5 + \sum_{x_i \in \mathbf{x}_M} \left((x_i - 0.5)^2 - \cos(20\pi(x_i - 0.5))\right)\right),
\end{eqnarray}
subject to $0 \le x_i \le 1$, for $i = 1,2,\dots,7$.\\
Here $\mathbf{x} = (x_1,x_2,x_3,x_4,x_5,x_6,x_7)$ and $\mathbf{x}_M = (x_3,x_4,x_5,x_6,x_7)$.

\item DTLZ2

Minimize $f_1$,$f_2$,$f_3$, where
\begin{eqnarray}
f_1(\mathbf{x}) &=& (1 + g(\mathbf{x}_M))\cos(x_1\pi/2)\cos(x_2\pi/2), \\
f_2(\mathbf{x}) &=& (1 + g(\mathbf{x}_M))\cos(x_1\pi/2)\sin(x_2\pi/2), \\
f_3(\mathbf{x}) &=& (1 + g(\mathbf{x}_M))\sin(x_1\pi/2), \\
g(\mathbf{x}_M) &=& \sum_{x_i \in \mathbf{x}_M} (x_i - 0.5)^2,
\end{eqnarray}
subject to $0 \le x_i \le 1$, for $i = 1,2,\dots,12$.\\
Here $\mathbf{x} = (x_1,x_2,\dots,x_{12})$ and $\mathbf{x}_M = (x_3,x_4,\dots,x_{12})$.

\item DTLZ4

Minimize $f_1$,$f_2$,$f_3$, where
\begin{eqnarray}
f_1(\mathbf{x}) &=& (1 + g(\mathbf{x}_M))\cos(x_1^{100}\pi/2)\cos(x_2^{100}\pi/2), \\
f_2(\mathbf{x}) &=& (1 + g(\mathbf{x}_M))\cos(x_1^{100}\pi/2)\sin(x_2^{100}\pi/2), \\
f_3(\mathbf{x}) &=& (1 + g(\mathbf{x}_M))\sin(x_1^{100}\pi/2), \\
g(\mathbf{x}_M) &=& \sum_{x_i \in \mathbf{x}_M} (x_i - 0.5)^2,
\end{eqnarray}
subject to $0 \le x_i \le 1$, for $i = 1,2,\dots,12$.\\
Here $\mathbf{x} = (x_1,x_2,\dots,x_{12})$ and $\mathbf{x}_M = (x_3,x_4,\dots,x_{12})$.

\end{itemize}

\begin{figure}[p]
\centering
\includegraphics[width=65mm]{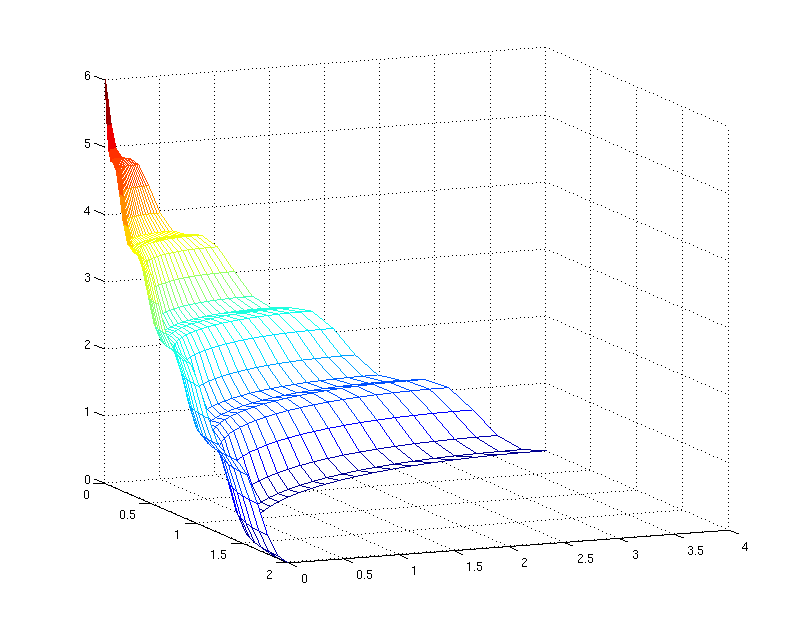}
\includegraphics[width=65mm]{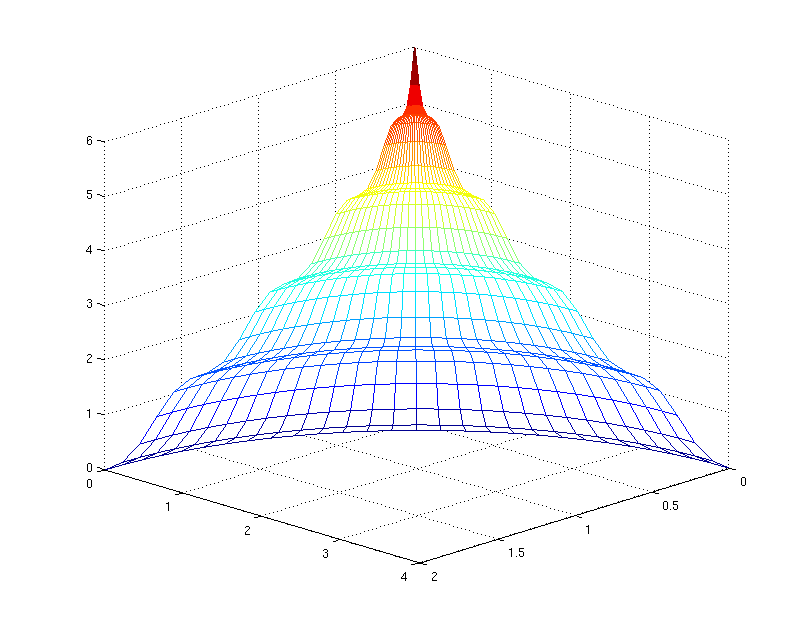}
\caption{\label{fig:WFG1_front} Exact Pareto front for problem WFG1.}
\end{figure}

Problem WFG1 is a benchmark problems introduced in \cite{Huband-2005-SMO}.
The following form is considered:
\begin{itemize}

\item WFG1

Minimize $f_1$,$f_2$,$f_3$, where
\begin{eqnarray}
f_1(\mathbf{x}) &=& 2\left[ (1-\cos(z_1\pi/2))(1-\cos(z_2\pi/2))\right], \\
f_2(\mathbf{x}) &=& 4\left[ (1-\cos(z_1\pi/2))(1-\sin(z_2\pi/2))\right], \\
f_3(\mathbf{x}) &=& 6\left[ 1 - z_1 - \frac{\cos(10\pi z_1 + \pi/2)}{10\pi}\right],
\end{eqnarray}
where $z_1 = z_1(\mathbf{x})$, and $z_2 = z_2(\mathbf{x})$
are auxiliary variables obtained from design variables $\mathbf{x} = (x_1,x_2,\dots,x_{24})$ 
by a series of nonlinear transformations (see \cite{Huband-2005-SMO} for their definition).
However, there is a slight difference in our application of transformation \texttt{b\_poly},
which we use with exponent 0.2 instead of 0.02 suggested in \cite{Huband-2005-SMO} due to 
numerical issues in floating point arithmetic.
Design variables have range $0 \le x_i \le 2i$, for $i = 1,2,\dots,24$.

The exact front of this problem is visualized in Figure \ref{fig:WFG1_front}.

\end{itemize}

\subsection{Metrics}
\label{sec:metrics}

The results are evaluated according to three measures.
The distance of members of the Pareto archive to the true Pareto front is measured using 
the \textbf{generational distance (GD)} \cite{Veldhuizen-1998-MEA}, which is defined as
\begin{equation}
\mbox{GD} = \sqrt{\frac{1}{n}\sum_{i=1}^{n}d_i^2},
\end{equation}
where $n$ is the number of nondominated solutions found by an~algorithm,
and $d_i$ is the Euclidean distance of the $i$-th solution to the exact front.
In order to evaluate the distance accurately, 
we implemented an~approach, 
that is able to iteratively find the closest point of the exact front for each approximate solution,
provided the analytic expression of the exact front is known.
This point is then used for measuring the distance.
Zero value of GD corresponds to all members of the archive on the exact front.

We evaluate also another measure of convergence,
denoted as \textbf{TOL5}.
It is defined as the lowest value, such that $d_i > \mbox{TOL5}$ holds for at most 5 \% of individuals.
In statistics, it is called the 95-th percentile with respect to distance.
Again, the lower value of TOL5 the better convergence.
Zero value indicates, that at least 95\% of archive members are on the exact front.
This metric is less sensitive to remote individuals than the GD value.

The uniformity of distribution of archive members is measured by \textbf{spacing} defined in \cite{Coello-2004-HMO}.
It is based on the distance to the nearest neighbour for each member of the archive, which is defined as
\begin{equation}
dn_i = \min_{\substack{j \in P \\ j\ne i}}\objdist{\objf_i}{\objf_j}.
\end{equation}
Now spacing is the ratio of standard deviation of values of these squared distances and their average,
i.e.
\begin{equation}
\mbox{spacing} = \frac{1}{\overline{dn}}\sqrt{\frac{1}{n-1}\sum_{i=1}^n(dn_i - \overline{dn})^2},
\end{equation}
where $\overline{dn}$ stands for the mean value 
\begin{equation}
\overline{dn} = \frac{1}{n}\sum_{i=1}^ndn_i .
\end{equation}
Consequently, zero spacing corresponds to uniform distribution of distances to the nearest neighbour.
Although this does not assure global uniformity of distribution (e.g. for pairs of individuals),
our experience with this metric is satisfactory.

The coverage of the Pareto front is not evaluated by means of a~metric,
but rather compared qualitatively at presented plots.

\subsection{Setting of algorithms}
\label{sec:setting}

All the algorithms from PISA package \cite{Bleuler-2003-PIS} (in the PISA context called selectors), 
i.e. NSGA-II, SPEA2, and IBEA, are used with the following setting of variator DTLZ:
\begin{itemize}
\item individual mutation probability \dots 1,
\item individual recombination probability \dots 1,
\item variable mutation probability \dots 0.01 ,
\item variable swap probability \dots 0.5,
\item variable recombination probability \dots 1,
\item $\eta$ mutation \dots 20,
\item $\eta$ recombination \dots 15,
\item use symmetric recombination \dots 1,
\end{itemize}
For variator WFG, these values are the same except the value
of variable mutation probability preset to 1 (default).
 
The simulations with PISA are performed with population of 100 individuals.
All of them are selected for mating,
producing 100 new individuals in each generation.
The tournament of 2 individuals is used in these selectors.
Experiments with IBEA are performed using the additive $\epsilon$-indicator with scaling factor $\kappa$ equal to 0.05.

The $\mu$ARMOGA is run with 4, 10 and 20 individuals in population, marked as $\mu$ARMOGA(4), $\mu$ARMOGA(10) and $\mu$ARMOGA(20), respectively
(for WFG1, only 4 members of population are considered).
The archive size is always set to 100 to produce results comparable with those of PISA algorithms.
Simple one-point crossover scheme without mutation is used.
It was reported by Oyama \cite{Oyama-2000-WDU}, that this scheme derived for binary coded algorithms \cite{Goldberg-1989-GAS} works reasonably well also for real-domain.
The version with 4 members uses reinitialization in each generation (DTLZ1, DTLZ2, DTLZ4) or once in four generations (WFG1), 
while for larger populations,
the reinitialization is performed once per 3 generations for all problems.
After reinitialization,
several existing archive members are put into the new population.
Their number is 2, 4 and 6 for the population of 4, 10 and 20 members, respectively.
Random selection of individuals for mating is then performed with this modified population. 
The other important parameters of $\mu$ARMOGA are set to the following values:
\begin{itemize}
\item adaptation factor $\delta $ \dots 1.4,
\item minimal standard deviation $\sigma_{min}$ \dots 0.8 (DTLZ1, WFG1), 0.005 (DTLZ2, DTLZ4),
\item recombination probability \dots 1.
\end{itemize}
Larger value of $\sigma_{min}$ 
helps to attain the global Pareto optimal front of multimodal problems such as DTLZ1  
and leads to faster convergence also in the case of WFG1.
In general, its large values emphasizes global exploration of the design space while small values lead to refined search.

\subsection{Experiments}
\label{sec:experiments}

\begin{figure}[p]
\centering
\includegraphics[width=58mm]{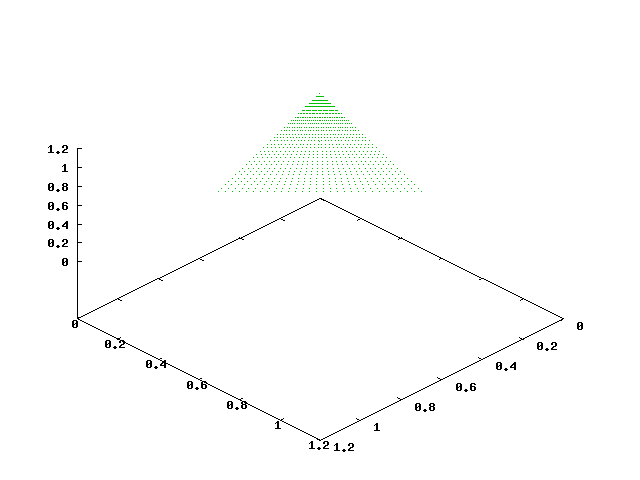}
\includegraphics[width=58mm]{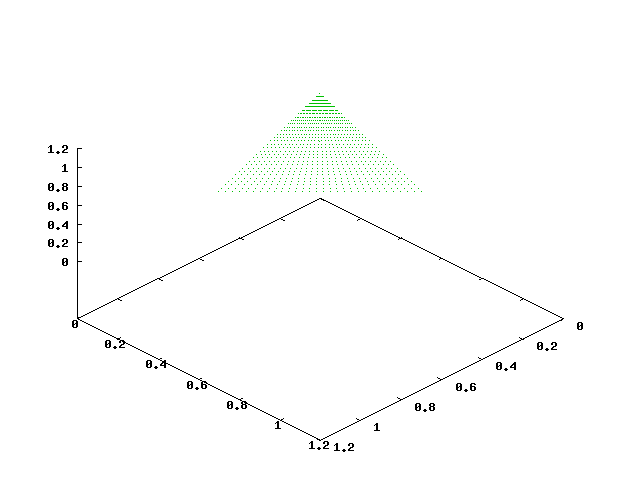} \\
\includegraphics[width=58mm]{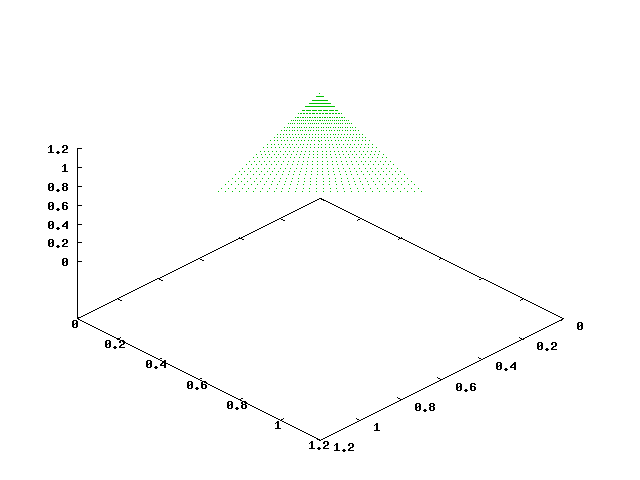} 
\includegraphics[width=58mm]{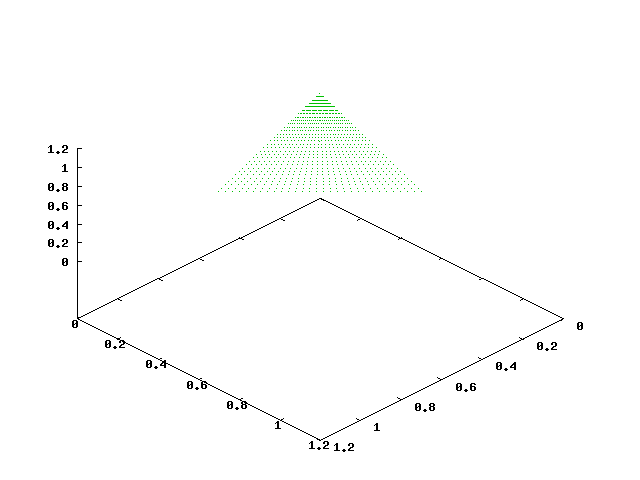} \\
\includegraphics[width=58mm]{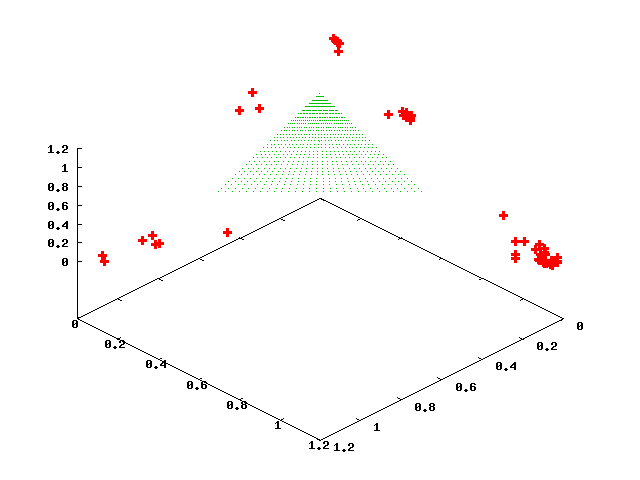}
\includegraphics[width=58mm]{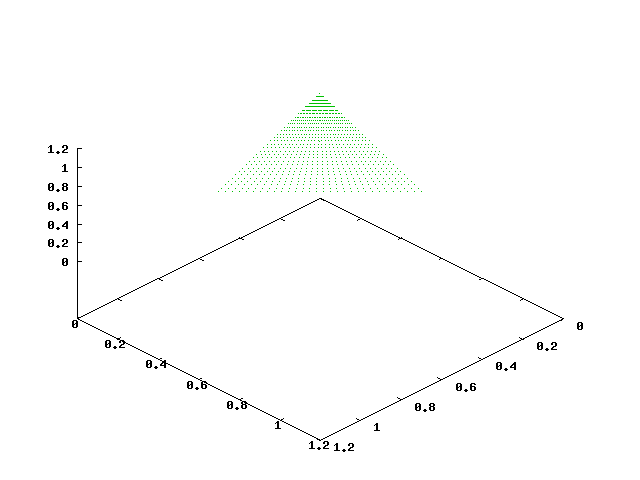}
\caption{\label{fig:DTLZ1_4000} Pareto front after 4\,000 function evaluations, problem DTLZ1, 
NSGA-II (top left), SPEA2 (centre left), IBEA (bottom left), and $\mu$ARMOGA with population size 4 (top right), 10 (centre right), 20 (bottom right).}
\end{figure}

\begin{figure}[p]
\centering
\includegraphics[width=58mm]{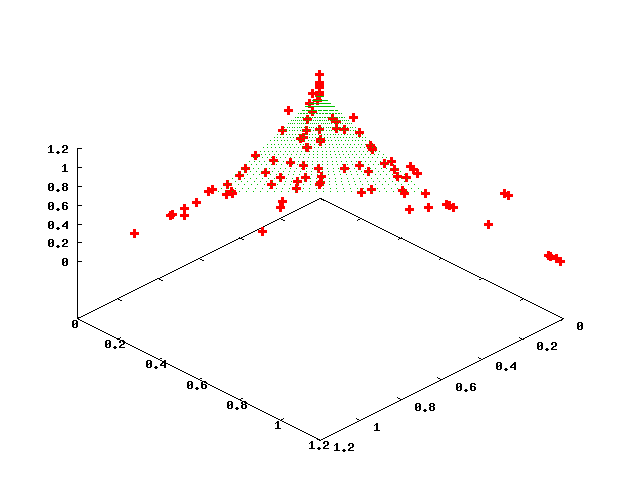}
\includegraphics[width=58mm]{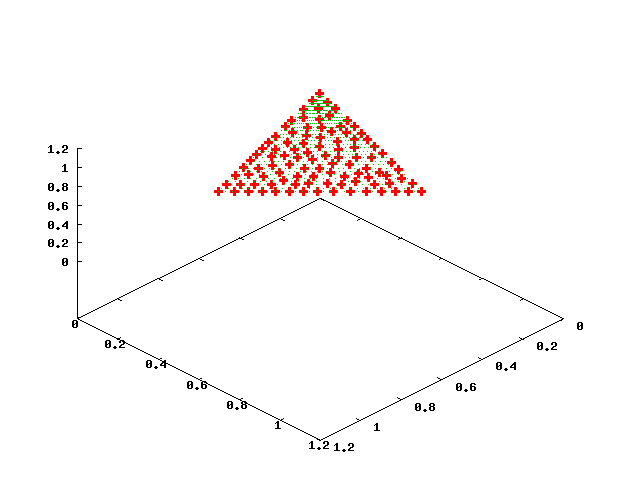} \\
\includegraphics[width=58mm]{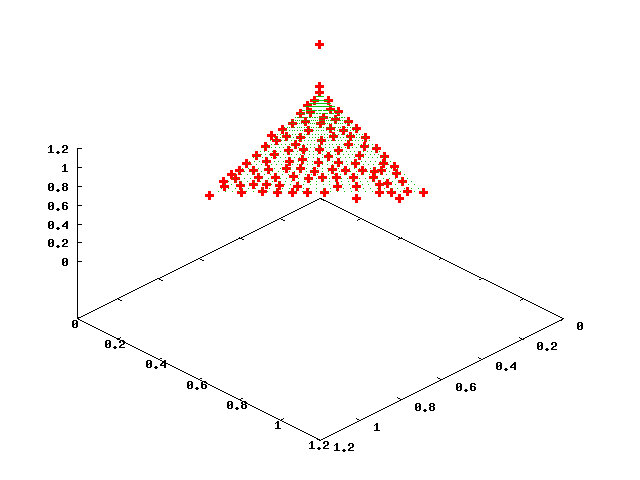} 
\includegraphics[width=58mm]{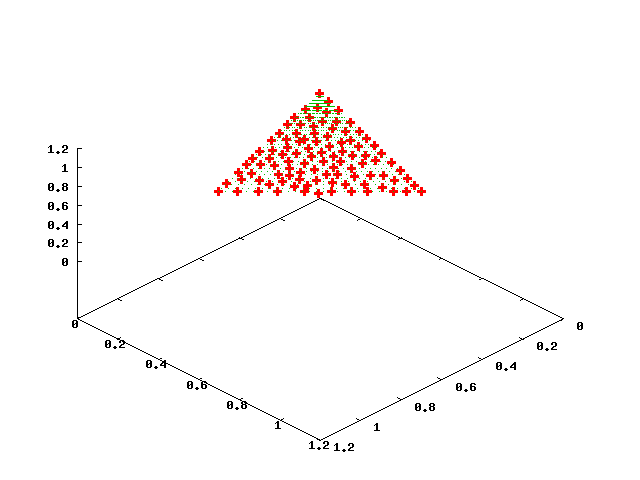} \\
\includegraphics[width=58mm]{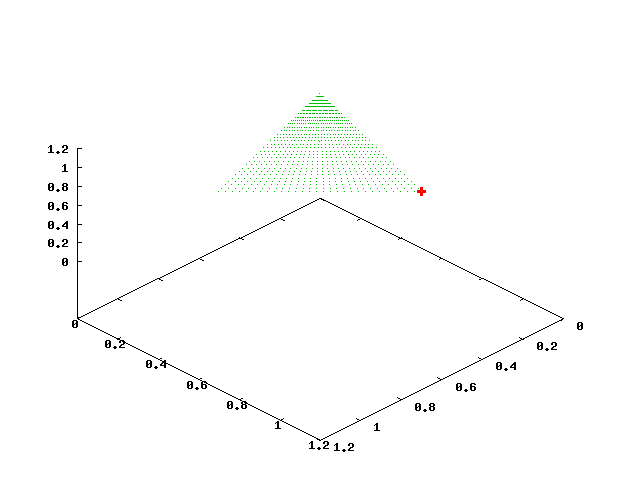}
\includegraphics[width=58mm]{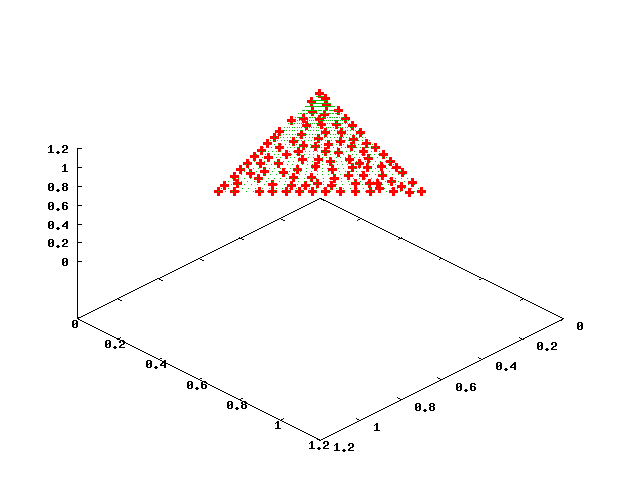}
\caption{\label{fig:DTLZ1_100000} Pareto front after 100\,000 function evaluations, problem DTLZ1,
NSGA-II (top left), SPEA2 (centre left), IBEA (bottom left), and $\mu$ARMOGA with population size 4 (top right), 10 (centre right), 20 (bottom right).}
\end{figure}

\begin{figure}[p]
\centering
\includegraphics[width=58mm]{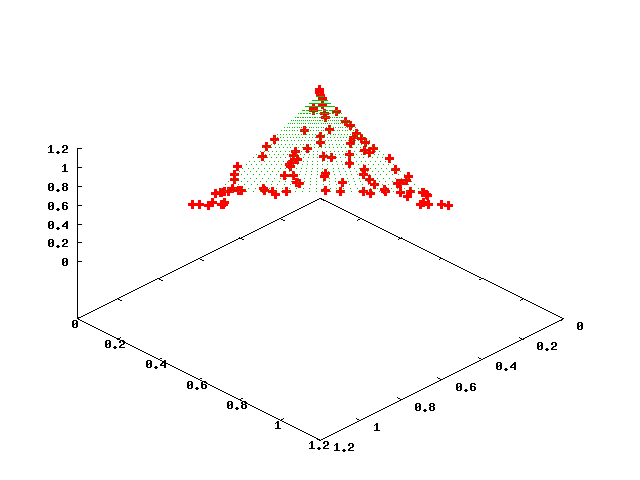}
\includegraphics[width=58mm]{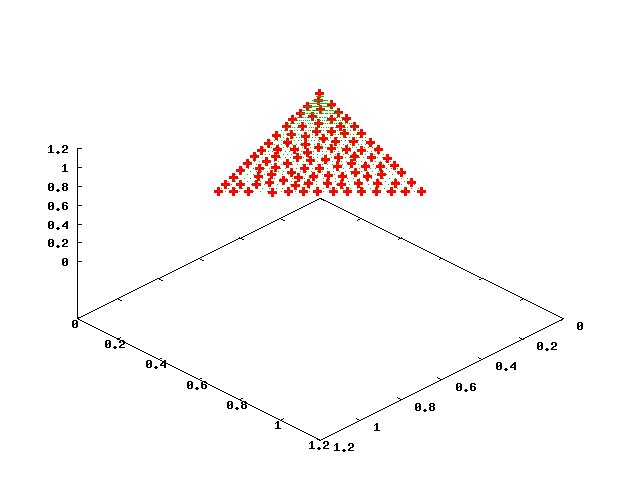} \\
\includegraphics[width=58mm]{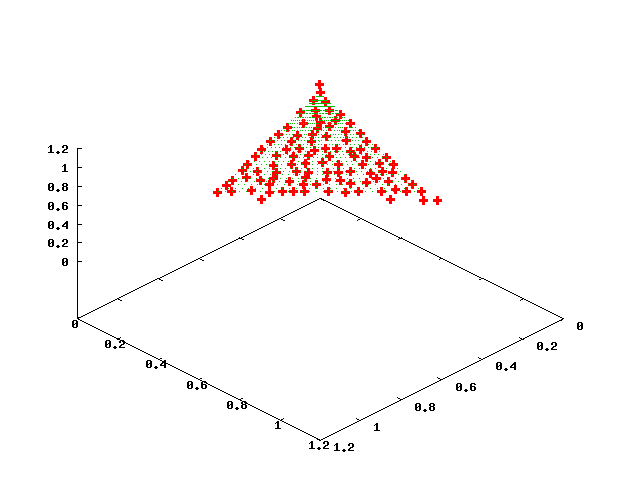} 
\includegraphics[width=58mm]{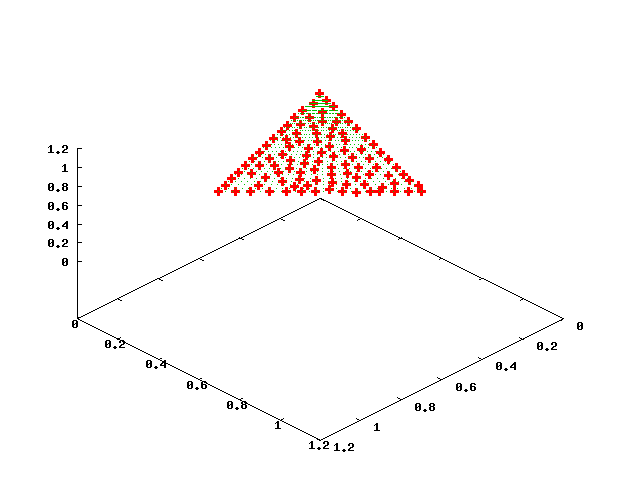} \\
\includegraphics[width=58mm]{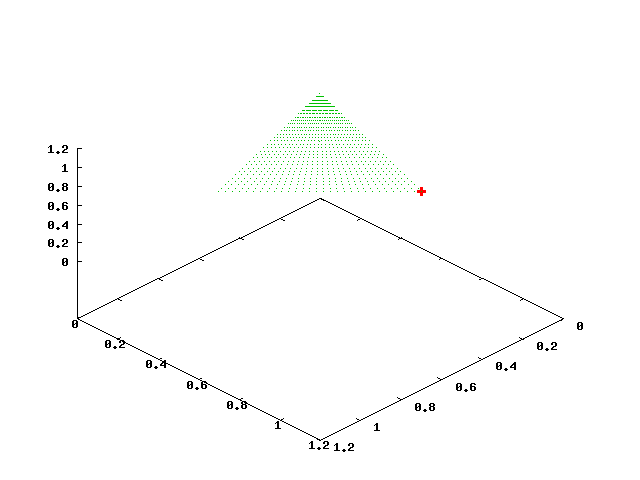}
\includegraphics[width=58mm]{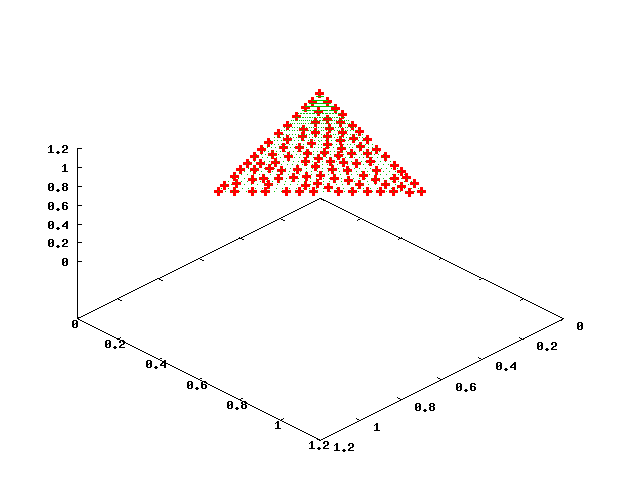}
\caption{\label{fig:DTLZ1_200000} Pareto front after 200\,000 function evaluations, problem DTLZ1,
NSGA-II (top left), SPEA2 (centre left), IBEA (bottom left), and $\mu$ARMOGA with population size 4 (top right), 10 (centre right), 20 (bottom right).}
\end{figure}

\begin{figure}[p]
\centering
\includegraphics[width=58mm]{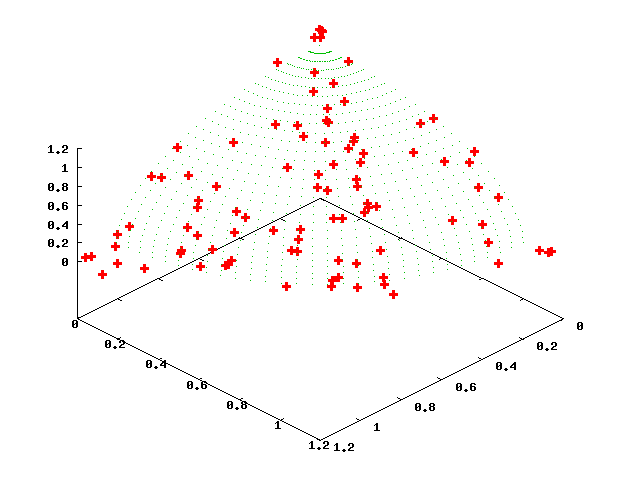}
\includegraphics[width=58mm]{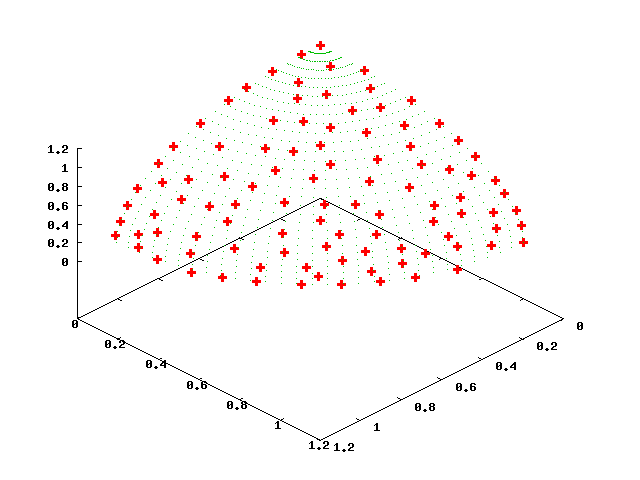} \\
\includegraphics[width=58mm]{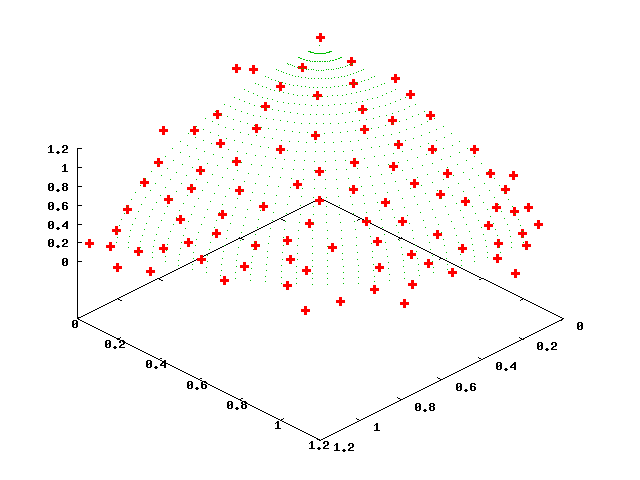} 
\includegraphics[width=58mm]{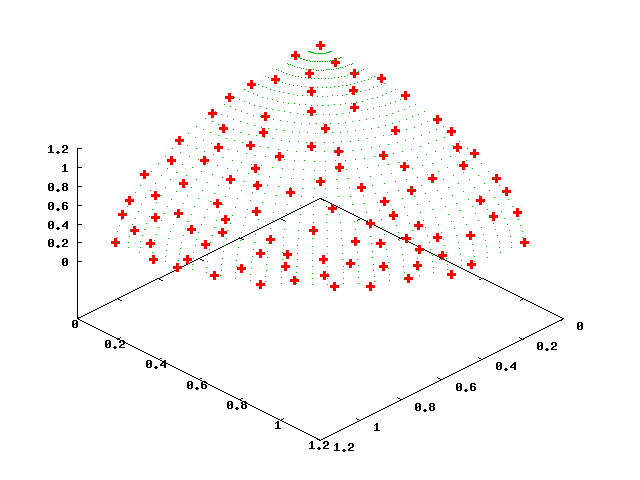} \\
\includegraphics[width=58mm]{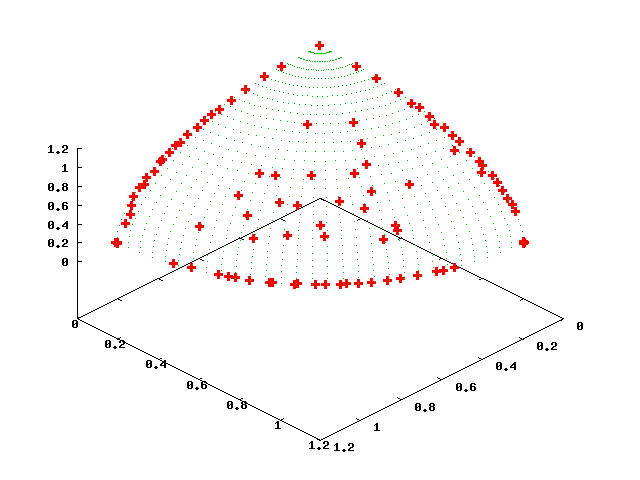}
\includegraphics[width=58mm]{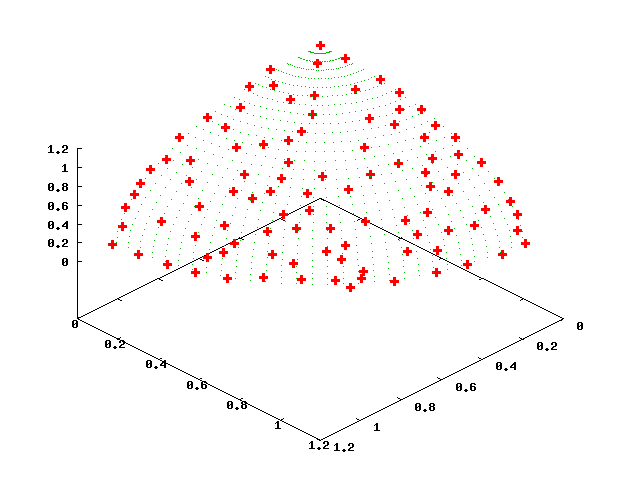}
\caption{\label{fig:DTLZ2_4000} Pareto front after 4\,000 function evaluations, problem DTLZ2, 
NSGA-II (top left), SPEA2 (centre left), IBEA (bottom left), and $\mu$ARMOGA with population size 4 (top right), 10 (centre right), 20 (bottom right).}
\end{figure}

\begin{figure}[p]
\centering
\includegraphics[width=58mm]{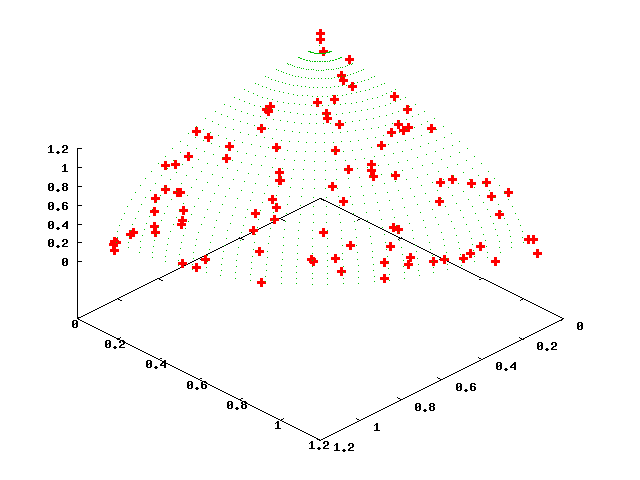}
\includegraphics[width=58mm]{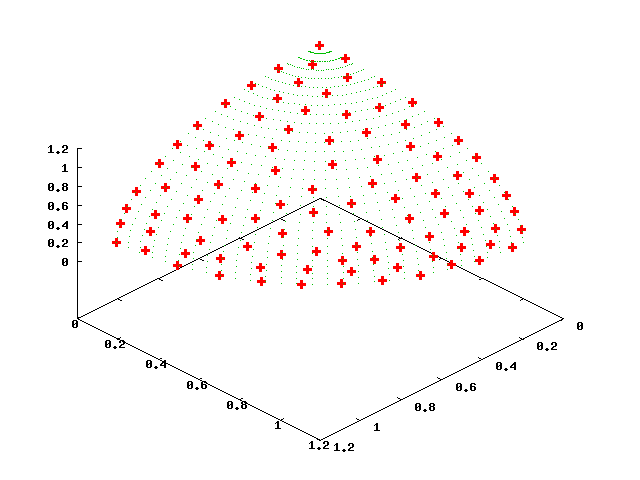} \\
\includegraphics[width=58mm]{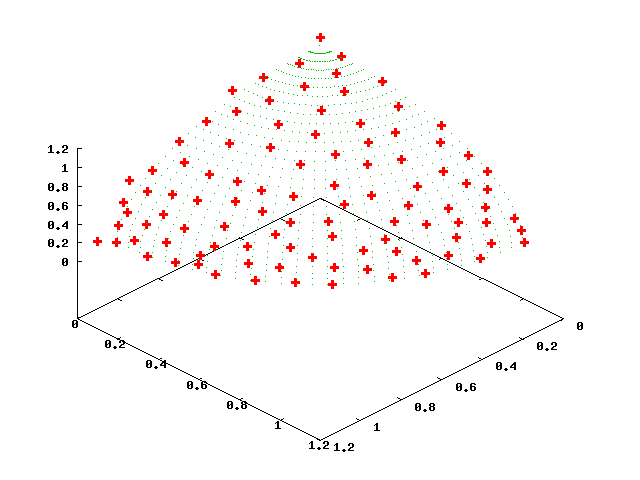} 
\includegraphics[width=58mm]{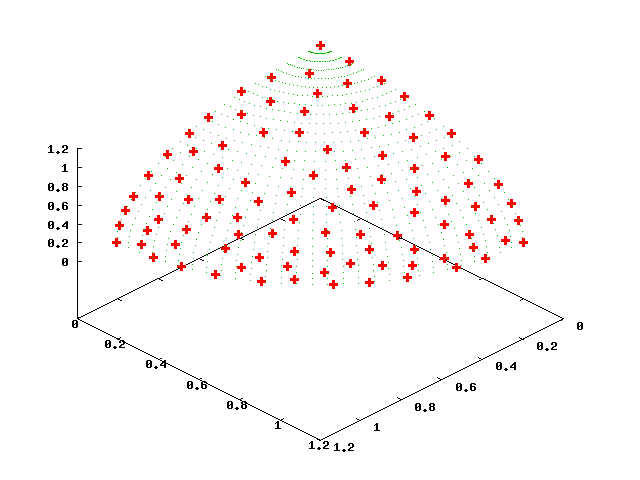} \\
\includegraphics[width=58mm]{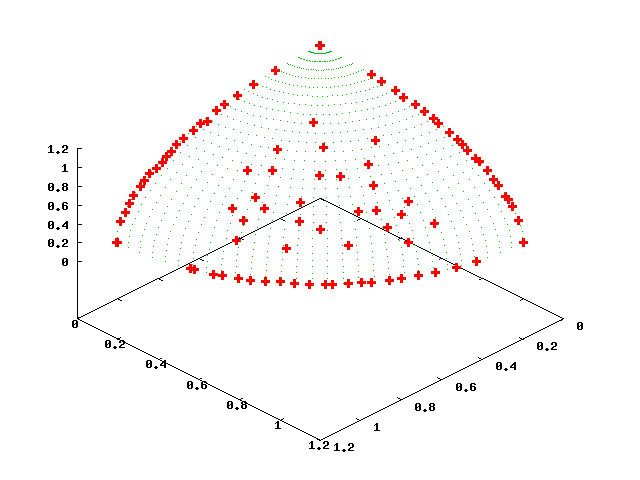}
\includegraphics[width=58mm]{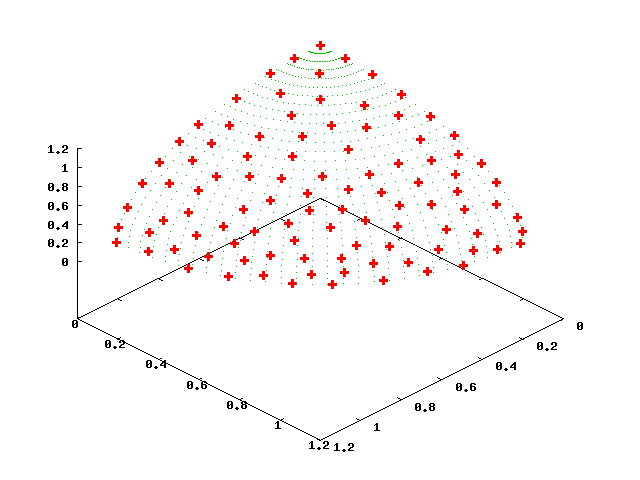}
\caption{\label{fig:DTLZ2_20000} Pareto front after 20\,000 function evaluations, problem DTLZ2, 
NSGA-II (top left), SPEA2 (centre left), IBEA (bottom left), and $\mu$ARMOGA with population size 4 (top right), 10 (centre right), 20 (bottom right).}
\end{figure}

\begin{figure}[p]
\centering
\includegraphics[width=58mm]{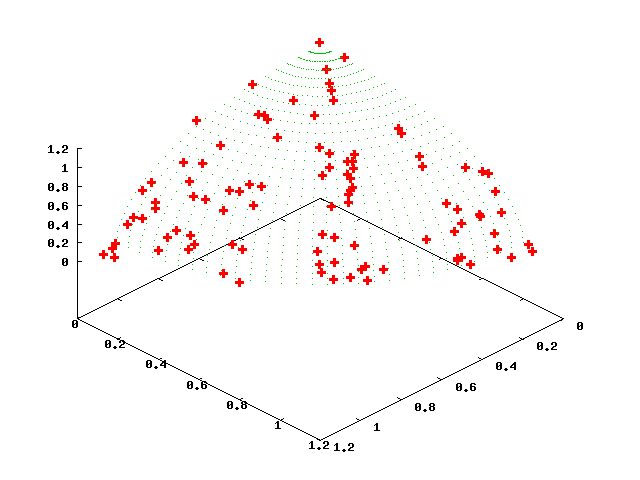}
\includegraphics[width=58mm]{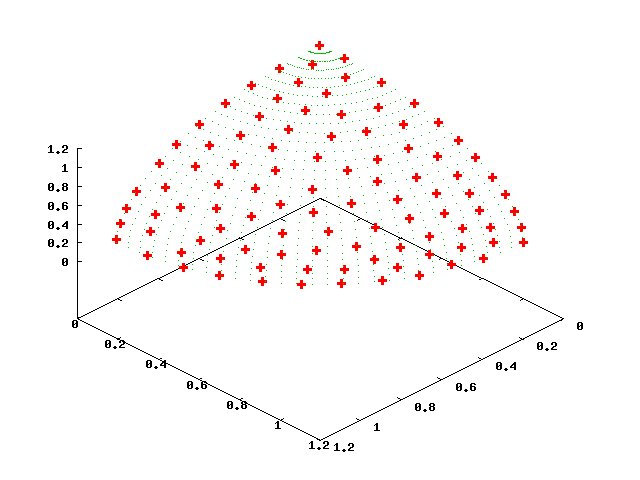} \\
\includegraphics[width=58mm]{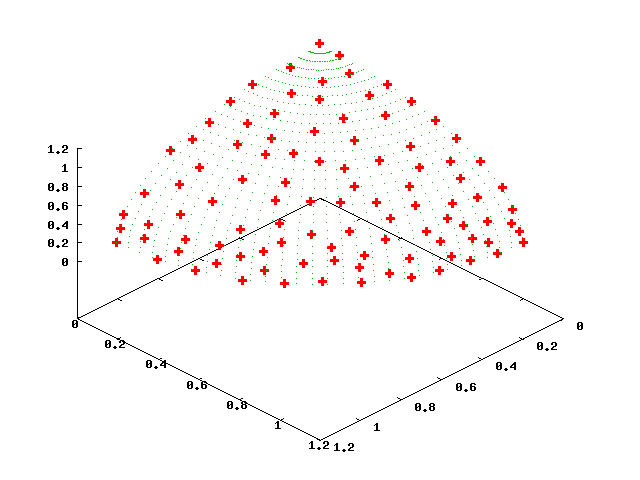} 
\includegraphics[width=58mm]{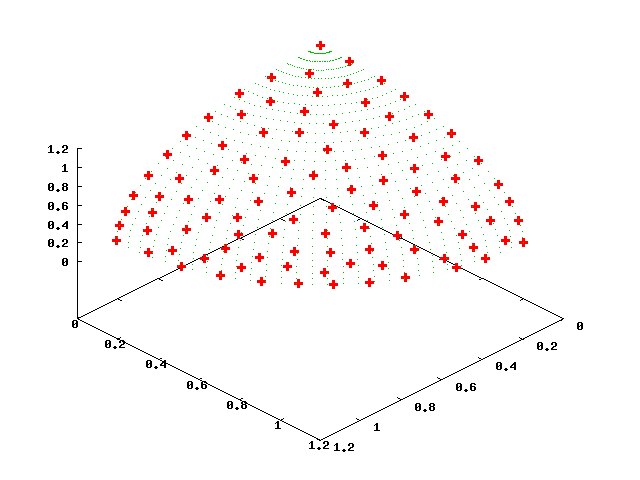} \\
\includegraphics[width=58mm]{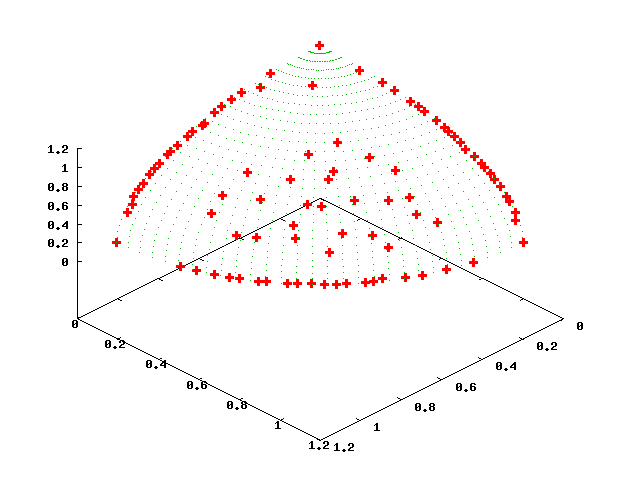}
\includegraphics[width=58mm]{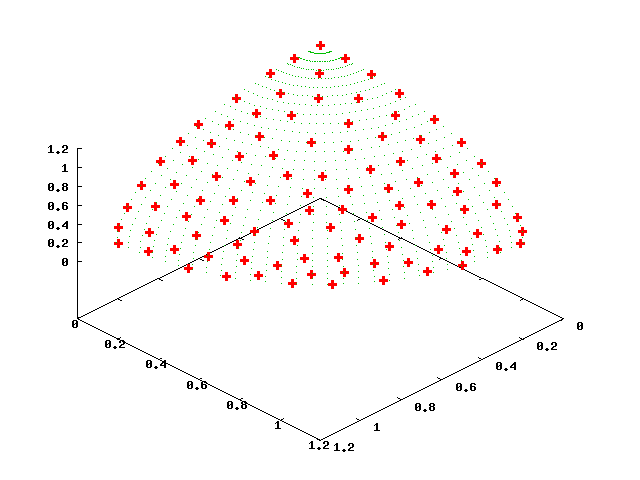}
\caption{\label{fig:DTLZ2_40000} Pareto front after 40\,000 function evaluations, problem DTLZ2, 
NSGA-II (top left), SPEA2 (centre left), IBEA (bottom left), and $\mu$ARMOGA with population size 4 (top right), 10 (centre right), 20 (bottom right).}
\end{figure}

\begin{figure}[p]
\centering
\includegraphics[width=58mm]{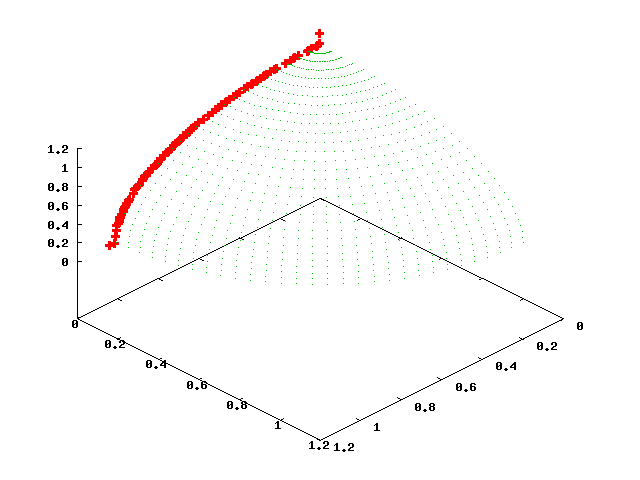}
\includegraphics[width=58mm]{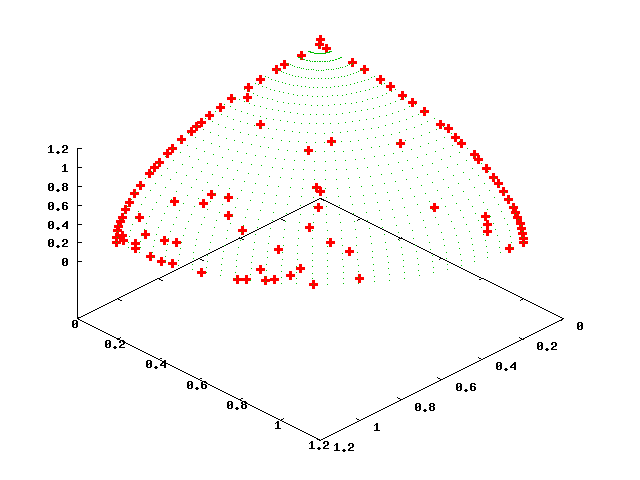} \\
\includegraphics[width=58mm]{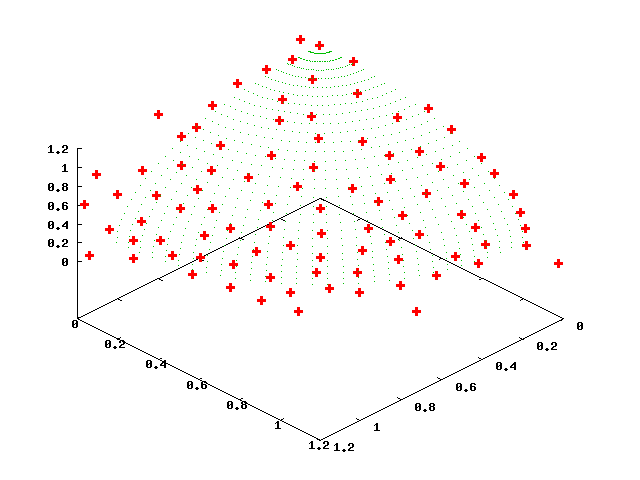} 
\includegraphics[width=58mm]{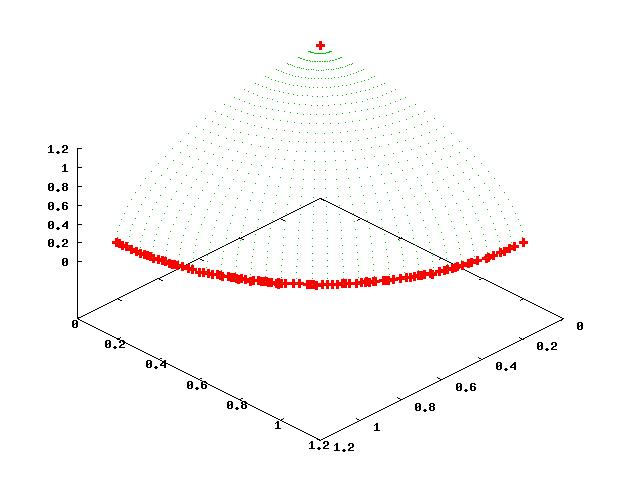} \\
\includegraphics[width=58mm]{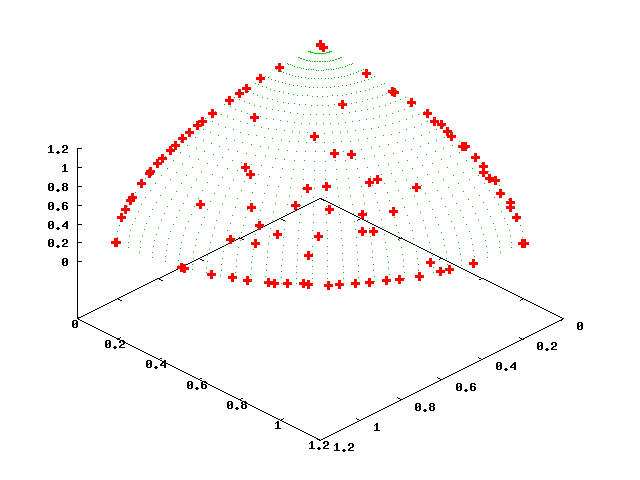}
\includegraphics[width=58mm]{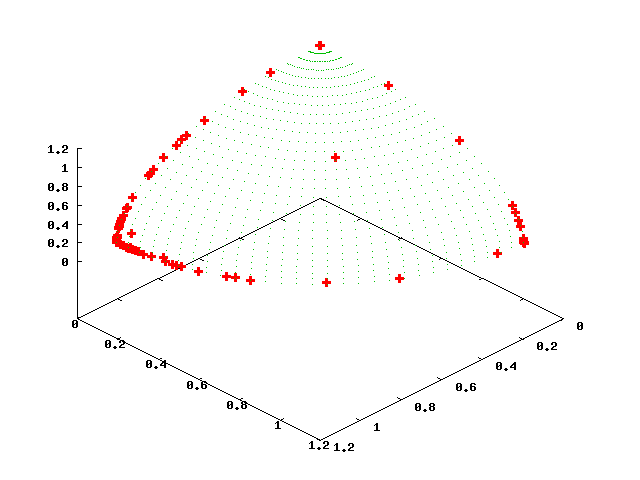}
\caption{\label{fig:DTLZ4_4000} Pareto front after 4\,000 function evaluations, problem DTLZ4, 
NSGA-II (top left), SPEA2 (centre left), IBEA (bottom left), and $\mu$ARMOGA with population size 4 (top right), 10 (centre right), 20 (bottom right).}
\end{figure}

\begin{figure}[p]
\centering
\includegraphics[width=58mm]{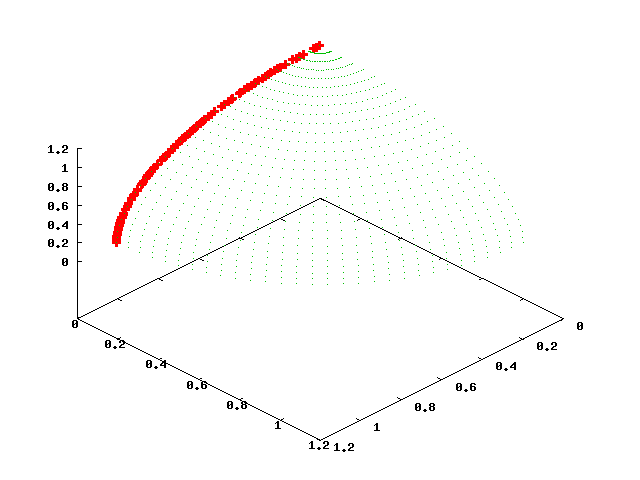}
\includegraphics[width=58mm]{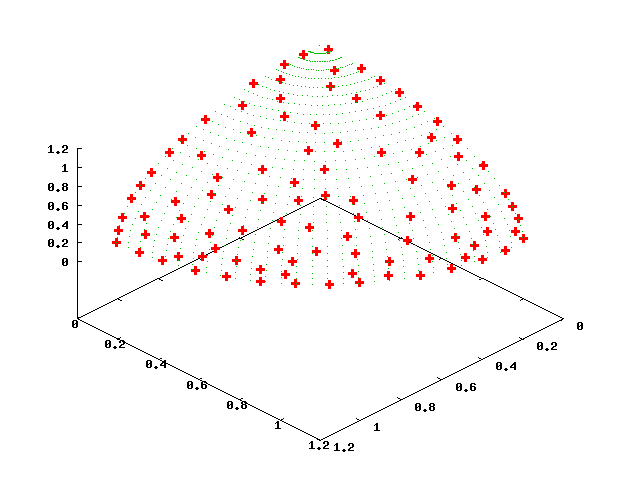} \\
\includegraphics[width=58mm]{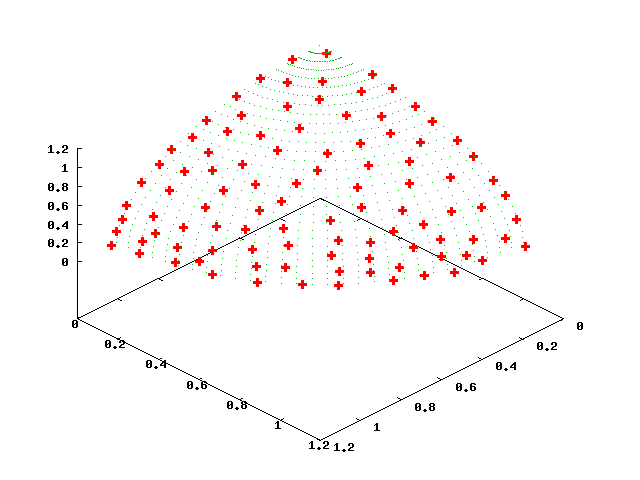} 
\includegraphics[width=58mm]{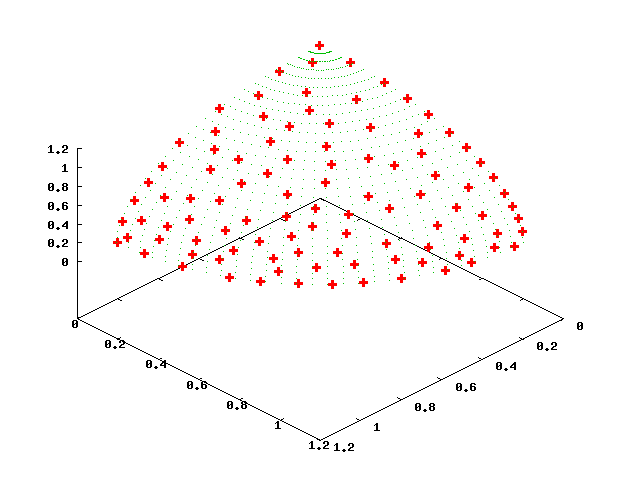} \\
\includegraphics[width=58mm]{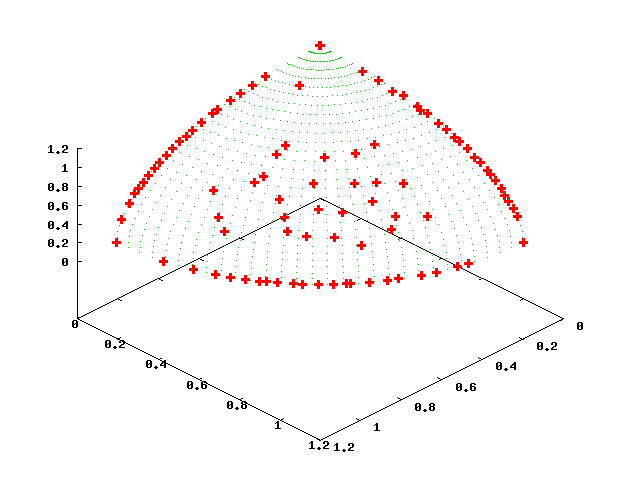}
\includegraphics[width=58mm]{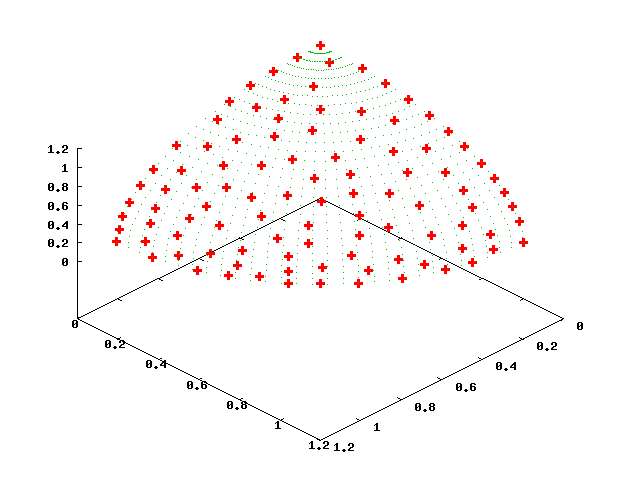}
\caption{\label{fig:DTLZ4_20000} Pareto front after 20\,000 function evaluations, problem DTLZ4, 
NSGA-II (top left), SPEA2 (centre left), IBEA (bottom left), and $\mu$ARMOGA with population size 4 (top right), 10 (centre right), 20 (bottom right).}
\end{figure}

\begin{figure}[p]
\centering
\includegraphics[width=58mm]{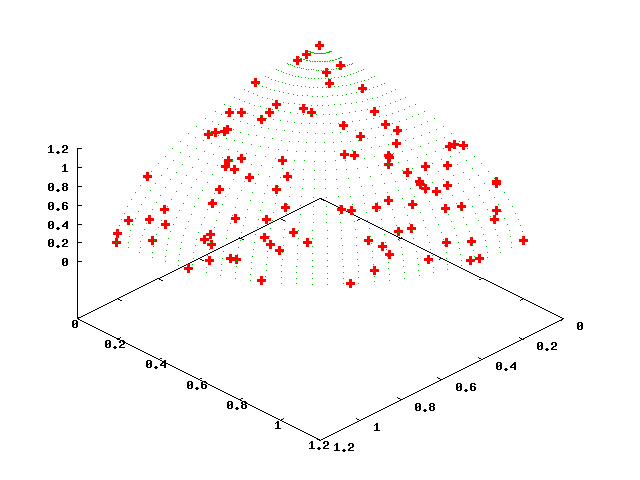}
\includegraphics[width=58mm]{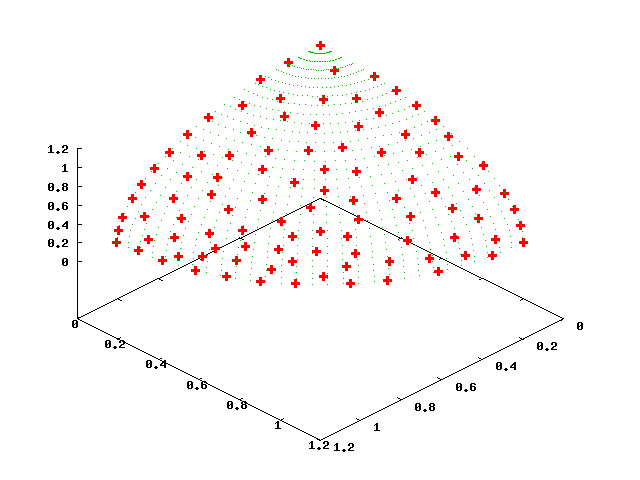} \\
\includegraphics[width=58mm]{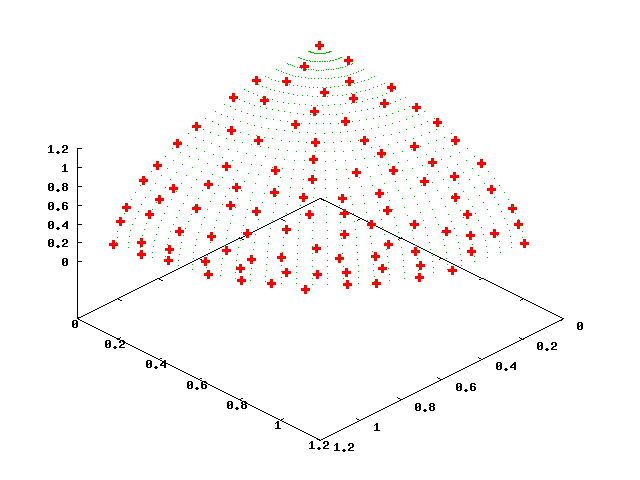} 
\includegraphics[width=58mm]{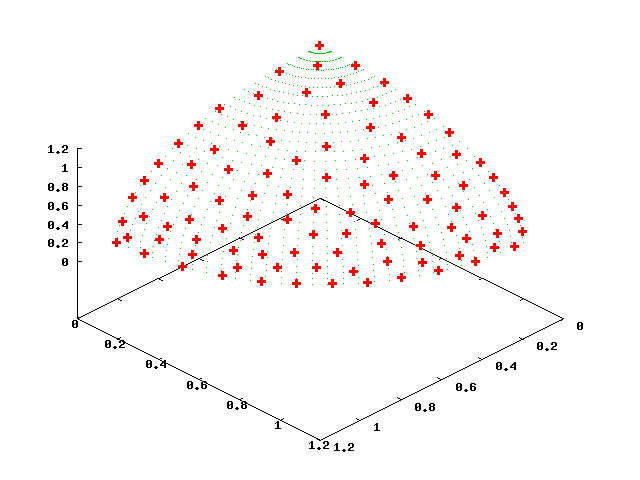} \\
\includegraphics[width=58mm]{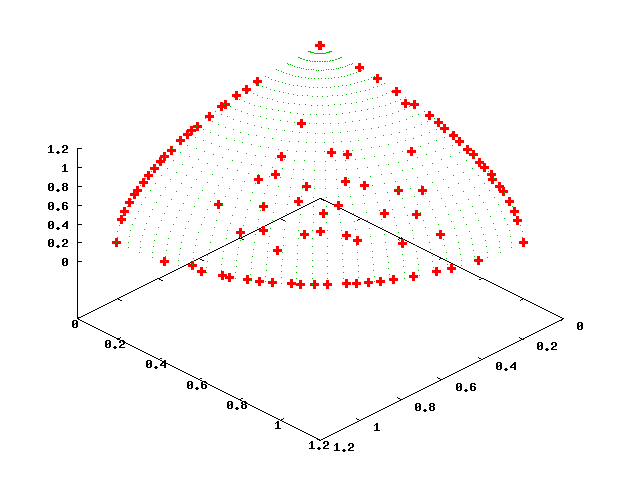}
\includegraphics[width=58mm]{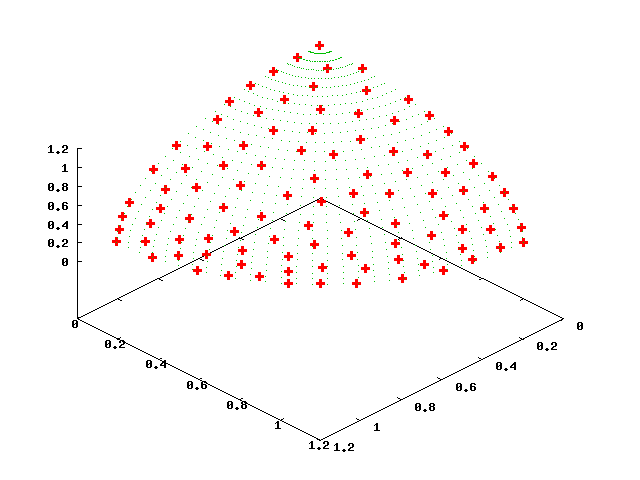}
\caption{\label{fig:DTLZ4_40000} Pareto front after 40\,000 function evaluations, problem DTLZ4, 
NSGA-II (top left), SPEA2 (centre left), IBEA (bottom left), and $\mu$ARMOGA with population size 4 (top right), 10 (centre right), 20 (bottom right).}
\end{figure}

\begin{figure}[p]
\centering
\includegraphics[width=58mm]{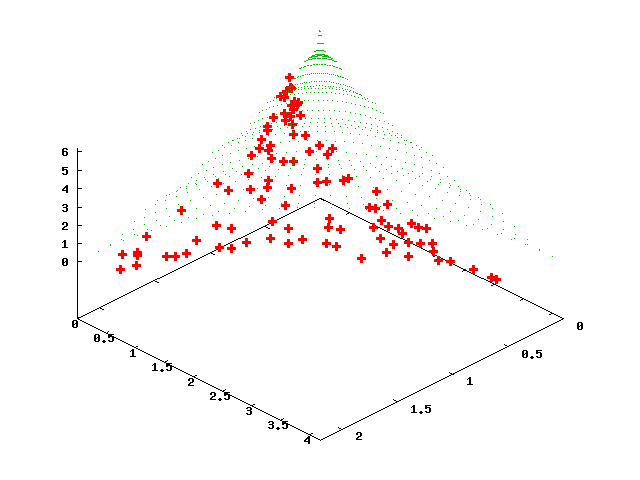}
\includegraphics[width=58mm]{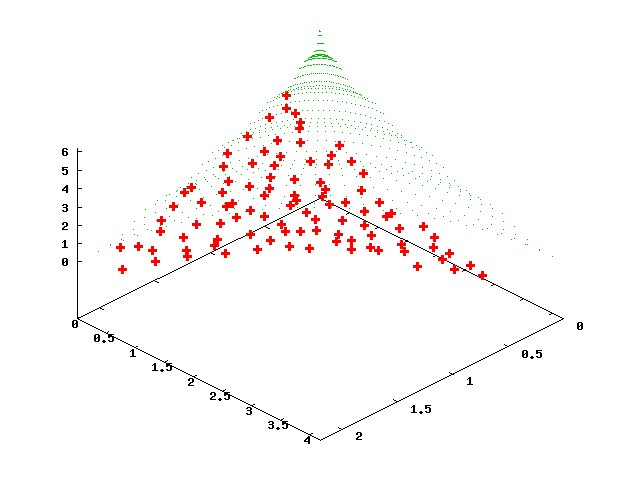} \\
\includegraphics[width=58mm]{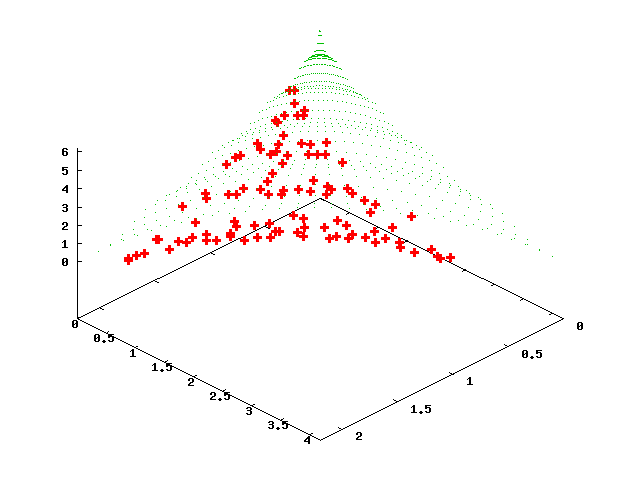}
\includegraphics[width=58mm]{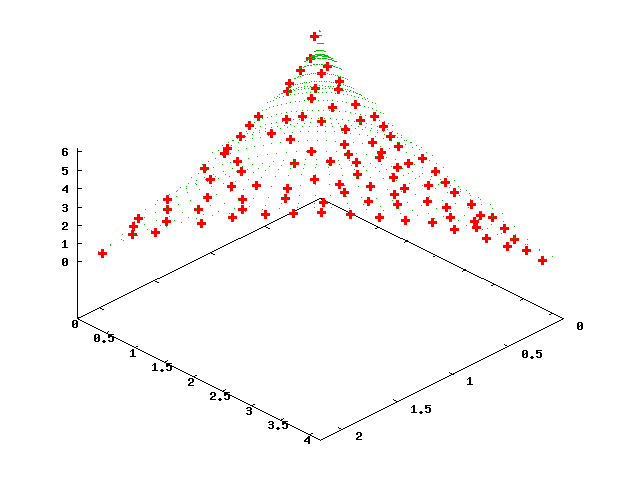}
\caption{\label{fig:WFG1_200000} Pareto front after 200\,000 function evaluations, problem WFG1, 
NSGA-II (top left), SPEA2 (top right), IBEA (bottom left), and $\mu$ARMOGA with population size 4 (bottom right).}
\end{figure}

\begin{figure}[p]
\centering
\includegraphics[width=58mm]{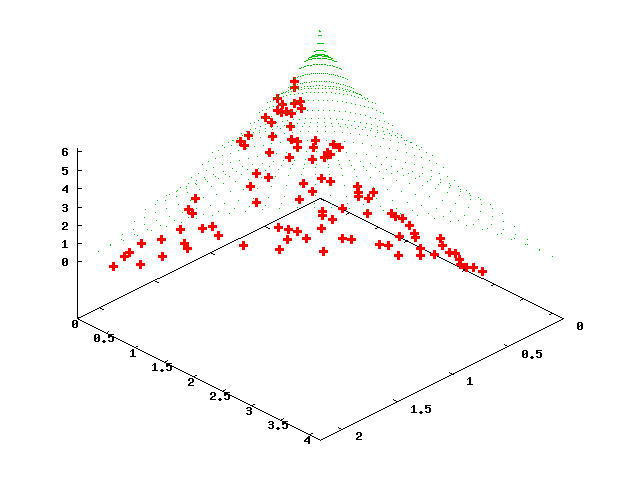}
\includegraphics[width=58mm]{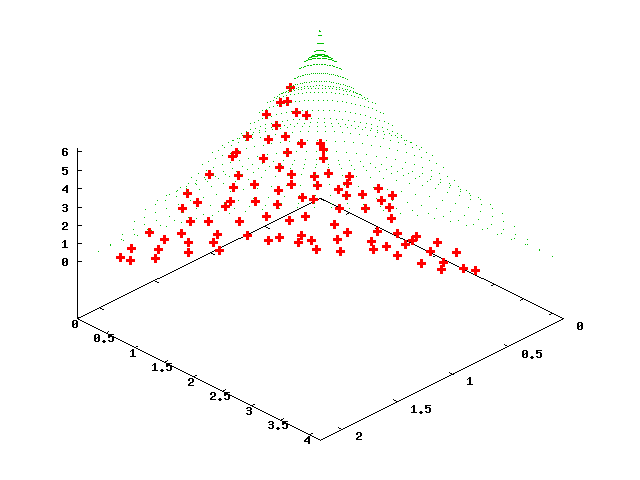} \\
\includegraphics[width=58mm]{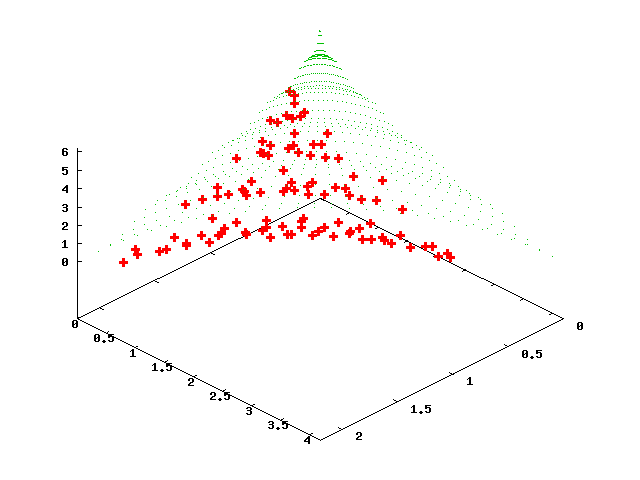}
\includegraphics[width=58mm]{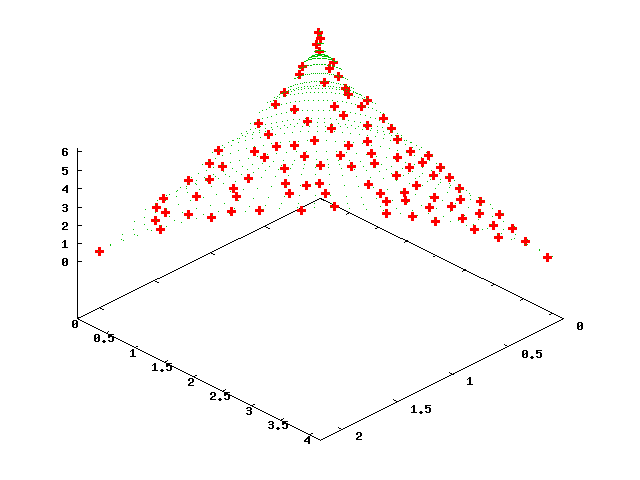}
\caption{\label{fig:WFG1_1000000} Pareto front after 1\,000\,000 function evaluations, problem WFG1, 
NSGA-II (top left), SPEA2 (top right), IBEA (bottom left), and $\mu$ARMOGA with population size 4 (bottom right).}
\end{figure}

\begin{figure}[p]
\centering
\includegraphics[width=58mm]{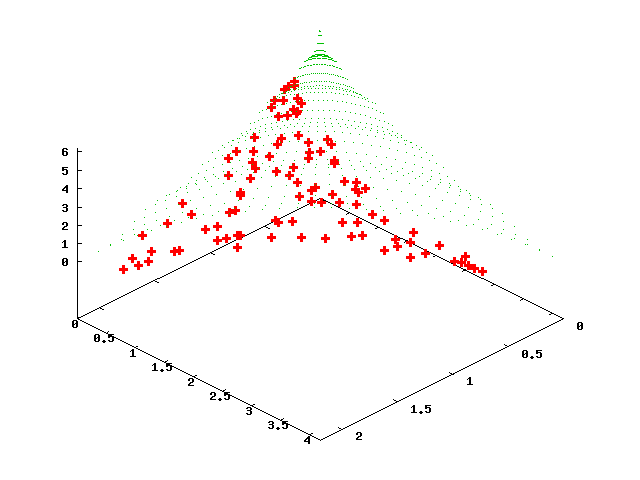}
\includegraphics[width=58mm]{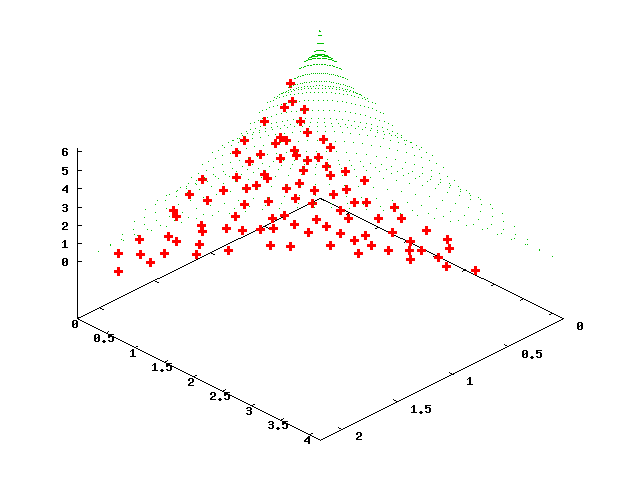} \\
\includegraphics[width=58mm]{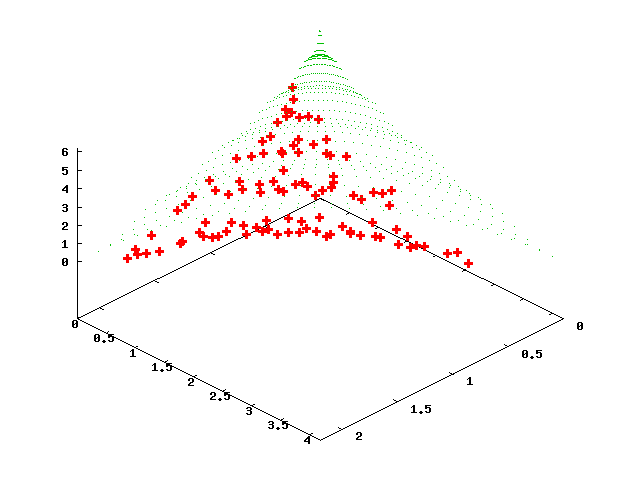}
\includegraphics[width=58mm]{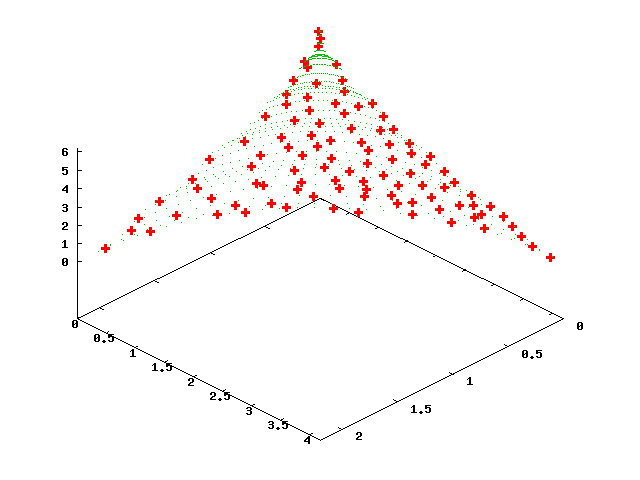}
\caption{\label{fig:WFG1_2000000} Pareto front after 2\,000\,000 function evaluations, problem WFG1, 
NSGA-II (top left), SPEA2 (top right), IBEA (bottom left), and $\mu$ARMOGA with population size 4 (bottom right).}
\end{figure}

The results for problems DTLZ1, DTLZ2 and DTLZ4 are summarized in Tables~\ref{tab:DTLZ1_4000}--\ref{tab:DTLZ4_40000},
and visualized in Figures~\ref{fig:DTLZ1_4000}--\ref{fig:DTLZ4_40000}.
For problem WFG1, results are summarized in Tables~\ref{tab:WFG1_4000}--\ref{tab:WFG1_2000000} 
and Figures~\ref{fig:WFG1_200000}--\ref{fig:WFG1_2000000}.
The values in tables are obtained as averages for 20 different seeds and where it makes sense, the best value is emphasized by bold font.
The approximation with the best distribution is selected out of the twenty runs of each algorithm for visualization.
The exact Pareto front is marked by grid of small dots in presented figures.
However, not all twenty seeds lead to a~successful approximation of the Pareto set in some instances,
especially for problem DTLZ4. 
Most of the algorithms suffer from problems with robustness with respect to initial population and produce \emph{degenerated} fronts for some seeds.
We consider a~front as degenerated, if all individuals have almost identical value of an~objective, 
and thus cover just a~line on the three dimensional surface of the exact front.
Numbers of degenerated fronts for all problems and methods are summarized in Table~\ref{tab:degenerated}. 

The efficiency of the proposed archiving technique is further demonstrated in comparison with the same approach, i.e. $\mu$ARMOGA, 
using crowding distance \cite{Deb-2001-MOO,Deb-2000-FEN}. 
That algorithm was described in \cite{Szollos-2008-AOM}.
Outputs of these experiments are summarized for DTLZ1 in Tables~\ref{tab:DTLZ1_macd} and \ref{tab:DTLZ1_mana}, 
for DTLZ2 in Tables~\ref{tab:DTLZ2_macd} and \ref{tab:DTLZ2_mana}, 
for DTLZ4 in Tables~\ref{tab:DTLZ4_macd} and \ref{tab:DTLZ4_mana}, 
and for WFG1 in Tables~\ref{tab:WFG1_macd} and \ref{tab:WFG1_mana}.
Obtained Pareto fronts are plotted in Figures~\ref{fig:DTLZ1_macd_mana_4000}--\ref{fig:WFG1_macd_mana_2000000}.

\section{Discussion of results}
\label{sec:discussion}

\subsection{DTLZ1}
This problem with three objectives has a linear Pareto optimal front that intersects the axes of the objective space at value 0.5.
Apart of the exact front, there exist a number of other parallel planes corresponding to local Pareto fronts.
As these also attract an evolution, problem DTLZ1 tests the ability of a genetic algorithm to cope with multi-modality.

As can be seen in Figure~\ref{fig:DTLZ1_4000},
none of the algorithms is able to reach the global Pareto front in 4\,000 evaluations for any seed, 
and metrics in Table~\ref{tab:DTLZ1_4000} do not provide much valuable information.
However, we can remark that IBEA and $\mu$ARMOGA(4) provide one order better convergence than the other algorithms and
$\mu$ARMOGA for all sizes of population provides reasonable spacing.

\begin{figure}[p]
\centering
\includegraphics[width=58mm]{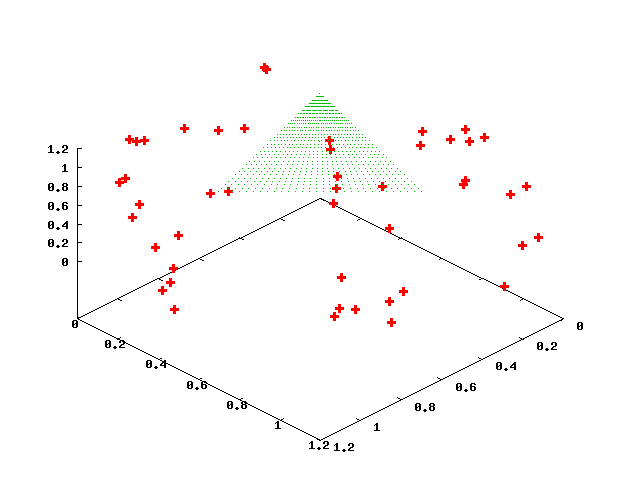}
\includegraphics[width=58mm]{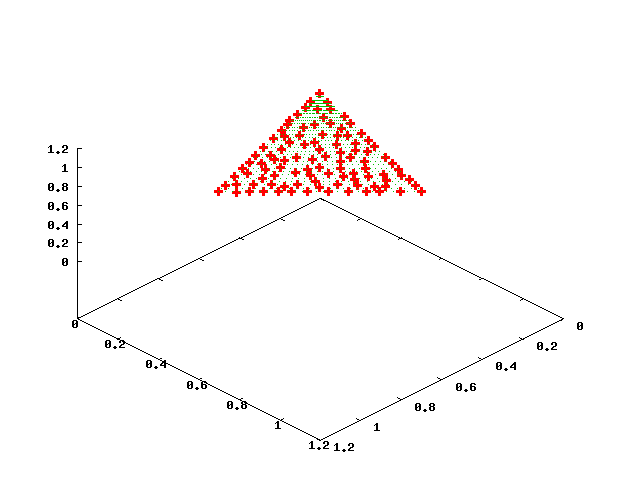} \\
\includegraphics[width=58mm]{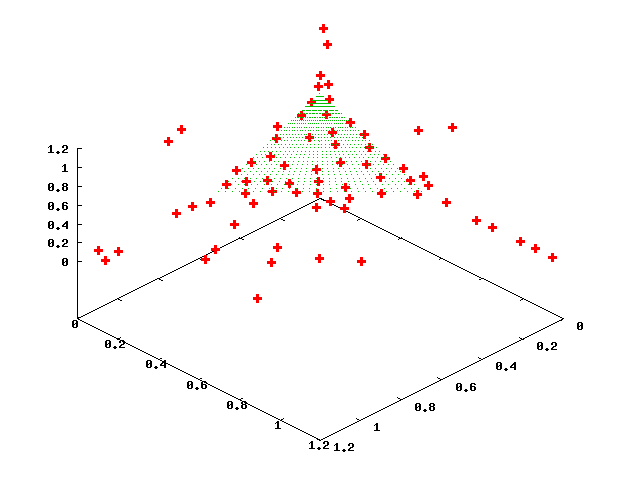} 
\includegraphics[width=58mm]{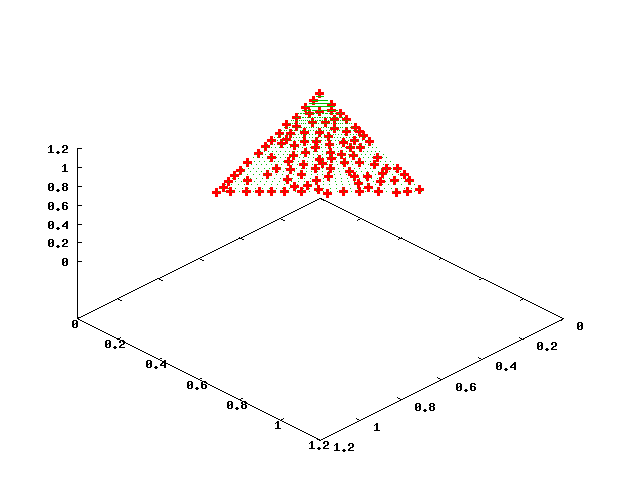} \\
\includegraphics[width=58mm]{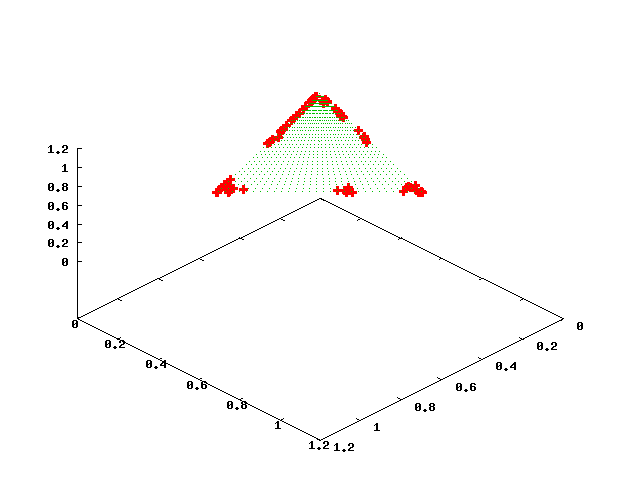}
\includegraphics[width=58mm]{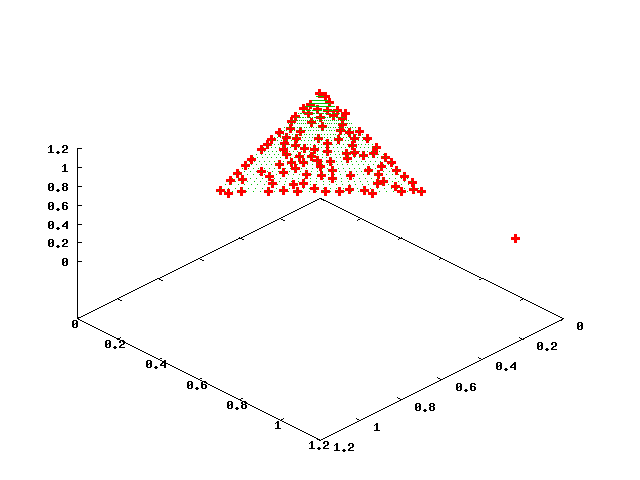}
\caption{\label{fig:DTLZ1_20000} Pareto front after 20\,000 function evaluations, problem DTLZ1, 
NSGA-II (top left), SPEA2 (centre left), IBEA (bottom left), and $\mu$ARMOGA with population size 4 (top right), 10 (centre right), 20 (bottom right).}
\end{figure}

However, all algorithms except NSGA-II are able to reach the global front for some seeds in 20\,000 evaluations (Figure~\ref{fig:DTLZ1_20000}).
Comparing visually the results in Figure~\ref{fig:DTLZ1_20000}, we can conclude that $\mu$ARMOGA(4) performs best,
which is supported by the best values of all three metrics in Table~\ref{tab:DTLZ1_20000}.
As individuals for many of the seeds are still away from the global front for all algorithms,
metrics in Table~\ref{tab:DTLZ1_20000} do not provide a detailed insight either.

\begin{figure}[p]
\centering
\includegraphics[width=58mm]{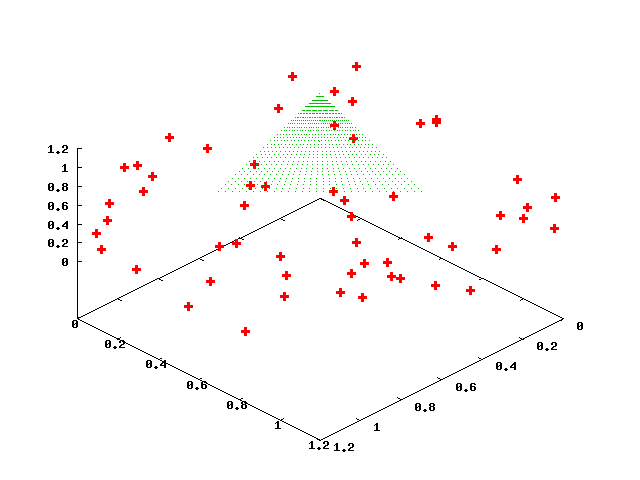}
\includegraphics[width=58mm]{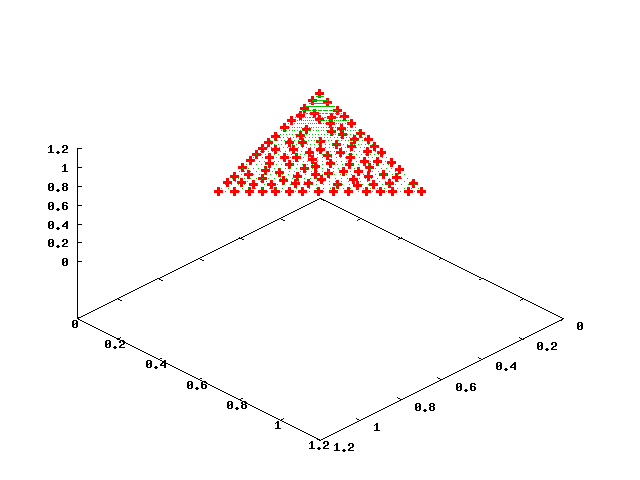} \\
\includegraphics[width=58mm]{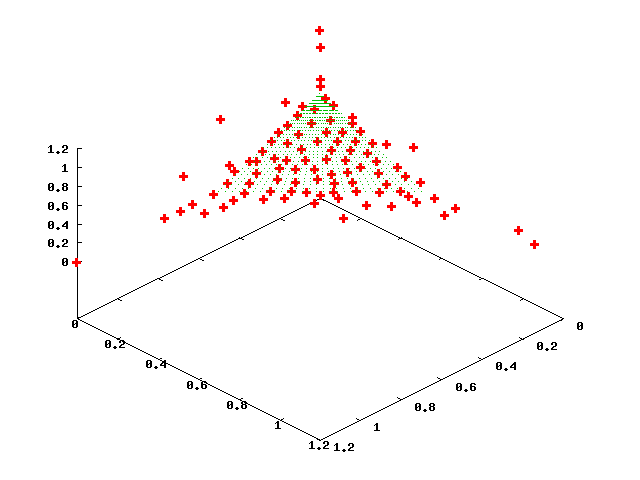} 
\includegraphics[width=58mm]{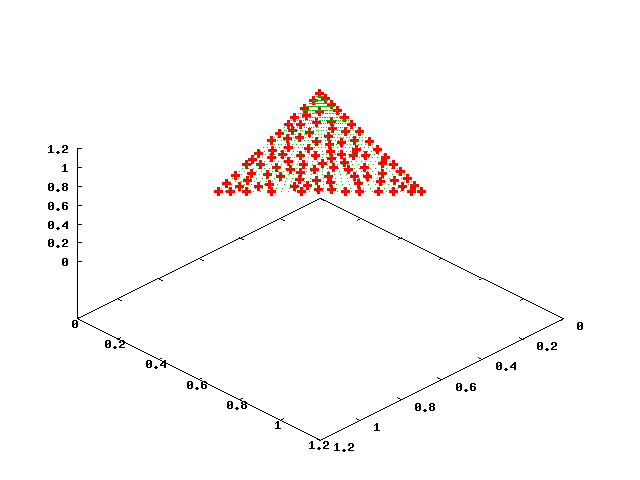} \\
\includegraphics[width=58mm]{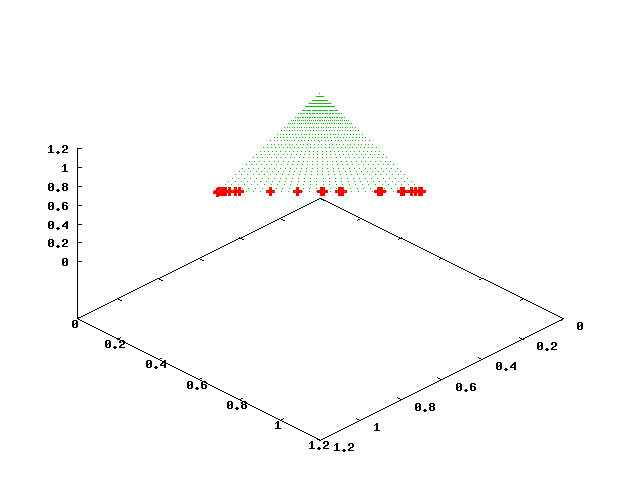}
\includegraphics[width=58mm]{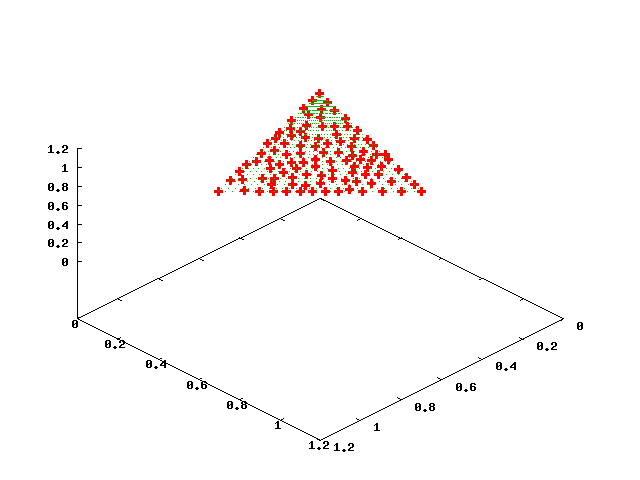}
\caption{\label{fig:DTLZ1_40000} Pareto front after 40\,000 function evaluations, problem DTLZ1,
NSGA-II (top left), SPEA2 (centre left), IBEA (bottom left), and $\mu$ARMOGA with population size 4 (top right), 10 (centre right), 20 (bottom right).}
\end{figure}

The situation is further improved with 40\,000 evaluations, 
for which all algorithms except NSGA-II are able to reach the global front for most of the seeds (Figure~\ref{fig:DTLZ1_40000}).
However, Figure~\ref{fig:DTLZ1_40000} shows that $\mu$ARMOGA (for all sizes of population) produces the best distribution,
which is confirmed by values of spacing in Table~\ref{tab:DTLZ1_40000}.
Since for some seeds the individuals still are not in vicinity of the true Pareto front, 
the averaged metrics in Table~\ref{tab:DTLZ1_40000} are still rather bad.
According to Table~\ref{tab:DTLZ1_40000}, the best convergence is in average attained by $\mu$ARMOGA(4) for this case.

For 100\,000 and 200\,000 evaluations, 
$\mu$ARMOGA(4) achieves the global Pareto-optimal front for all seeds.
All the other algorithms fail to find the global front for some seeds, 
which considerably spoils the metrics in Tables~\ref{tab:DTLZ1_100000} and \ref{tab:DTLZ1_200000}.
Since IBEA produces only degenerated fronts in these cases, 
metrics are not evaluated and are omitted in Tables~\ref{tab:DTLZ1_100000} and \ref{tab:DTLZ1_200000}.

Although the distribution of fronts obtained by $\mu$ARMOGA for all sizes of population is comparable to SPEA according to Figures~\ref{fig:DTLZ1_100000} and \ref{fig:DTLZ1_200000},
the metrics in Tables~\ref{tab:DTLZ1_100000} and \ref{tab:DTLZ1_200000} reveal that spacing is, in average, one order better by $\mu$ARMOGA than by SPEA.
The best average convergence metrics are obtained by $\mu$ARMOGA(4) (Tables~\ref{tab:DTLZ1_100000} and \ref{tab:DTLZ1_200000}).

\subsection{DTLZ2}
Problem DTLZ2 has three objectives, 
and the exact front corresponds to the part of a~unit sphere when restricted to the octant given by non-negative values of all three objectives.
This is the easiest problem for all compared algorithms and tests mainly the speed at which an algorithm is converging to the exact Pareto front.

Already for 4\,000 evaluations, the fronts obtained by all the compared methods are reasonably converged and distributed along the exact Pareto front.
Figure~\ref{fig:DTLZ2_4000} shows that $\mu$ARMOGA (regardless of the size of population) and SPEA produce the best distribution of individuals along the exact front,
whereas the distribution obtained by IBEA is rather poor.
This observation is confirmed by the spacing metric in Table~\ref{tab:DTLZ2_4000}.
The best convergence is achieved by IBEA according to GD and TOL5 metrics in Table~\ref{tab:DTLZ2_4000} followed by $\mu$ARMOGA.

Similar observations can be made from the results for 20\,000 evaluations (Table~\ref{tab:DTLZ2_20000} and Figure~\ref{fig:DTLZ2_20000})
and 40\,000 evaluations (Table~\ref{tab:DTLZ2_40000} and Figure~\ref{fig:DTLZ2_40000}) 
-- the best spacing is obtained for all sizes of population by $\mu$ARMOGA and the best convergence is attained by IBEA, 
although the distribution of individuals along the Pareto front is worse.

\subsection{DTLZ4}

Although the definition of problem DTLZ4 is similar to DTLZ2 (cf. Section \ref{sec:problems}),
the evolution is greatly influenced by the exponential transformation of design variables,
which maps most of the space towards the axes in design space.
This in turn pushes the evolution to the limits of the objective space.
Thus, problem DTLZ4 tests best of the three DTLZ problems the ability of a genetic algorithm to obtain uniform distribution of individuals along the Pareto optimal surface.

For 4\,000 evaluations, the best distribution is produced by SPEA2 (Figure~\ref{fig:DTLZ4_4000}). 
This is confirmed by results of spacing in Table~\ref{tab:DTLZ4_4000}.
However, $\mu$ARMOGA(4) produces the best converged results.

For 20\,000 evaluations, 
the distribution obtained by $\mu$ARMOGA is already visually comparable with SPEA2
in Figure~\ref{fig:DTLZ4_20000}.
Also spacing obtained by $\mu$ARMOGA is comparable to that of SPEA2 according to Table~\ref{tab:DTLZ4_20000} for 4 and 10 members of population.
Algorithms $\mu$ARMOGA(20), IBEA, and NSGA-II produce in average only slightly worse converged results.
The best GD and TOL5 values are attained by $\mu$ARMOGA(20).

In the case of 40\,000 evaluations, 
the best distribution of individuals is attained by $\mu$ARMOGA followed by SPEA2
according to Figure~\ref{fig:DTLZ4_40000} and also confirmed by values of spacing in Table~\ref{tab:DTLZ4_40000}.
The best convergence is again obtained by $\mu$ARMOGA(20), followed by $\mu$ARMOGA(10), $\mu$ARMOGA(4) and IBEA, respectively (Table~\ref{tab:DTLZ4_40000}).

\subsection{WFG1}

This is a difficult problem and all tested algorithms had problems with convergence to the Pareto front.
For this reason, number of evaluations of the objective function was increased to 2\,000\,000,
after which some algorithms were able to attain the exact front.

After 4\,000 evaluations, all algorithms produce results rather far from the Pareto optimal set (Table~\ref{tab:WFG1_4000}).
Nevertheles, SPEA2 produces the most uniform distribution according to the spacing metric.

After 20\,000 evaluations, $\mu$ARMOGA(4) slightly leads in convergence followed by IBEA (Table~\ref{tab:WFG1_20000}),
producing distribution with uniformity between SPEA2 (best spacing) and the rest of the algorithms. 
The same observations remain valid for 40\,000 evaluations (Table~\ref{tab:WFG1_40000}.

After 100\,000 as well as 200\,000 evaluations, $\mu$ARMOGA(4) dominates in convergence to the exact front 
(GD and TOL5 metrics in Tables~\ref{tab:WFG1_100000} and \ref{tab:WFG1_200000}),
producing distribution comparable with SPEA2 (best spacing).
However, Figure~\ref{fig:WFG1_200000} suggest, that $\mu$ARMOGA(4) covers the whole Pareto front, 
unlike SPEA2.

Letting the evolution run to 1\,000\,000 and 2\,000\,000 evaluations, $\mu$ARMOGA(4) dominates both in convergence 
(one order of magnitude compared to the second IBEA in Tables~\ref{tab:WFG1_100000} and \ref{tab:WFG1_200000})
and distribution along the exact Pareto front.
In spacing metric, $\mu$ARMOGA(4) is followed by SPEA2.
Figures~\ref{fig:WFG1_1000000} and \ref{fig:WFG1_2000000} show that distribution of individuals by $\mu$ARMOGA(4) uniformly covers 
the whole Pareto front, 
while the other algorithms approaches the region around $f_1 = 0$ only slowly.

\subsection{Comparison of crowding distance with the new archiving mechanism}

\begin{figure}[p]
\centering
\includegraphics[width=58mm]{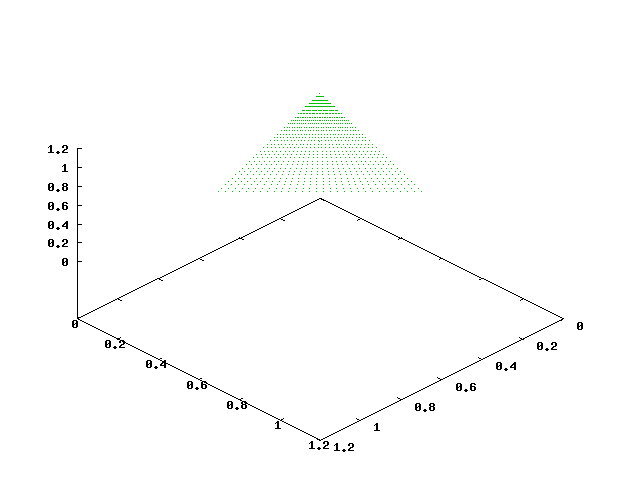}
\includegraphics[width=58mm]{emiarmoga4_DTLZ1_e4000} 
\caption{\label{fig:DTLZ1_macd_mana_4000} Pareto front after 4\,000 function evaluations, problem DTLZ1, 
$\mu$ARMOGA with population size 4 with crowding distance (left), and with the new proposed algorithm (right).}
\end{figure}

\begin{figure}[p]
\centering
\includegraphics[width=58mm]{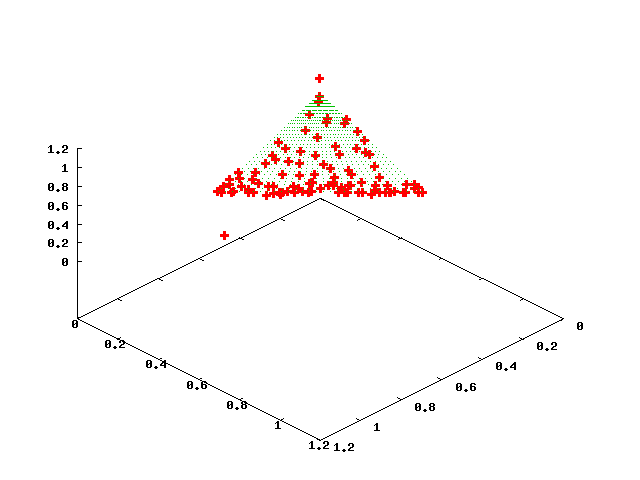}
\includegraphics[width=58mm]{emiarmoga4_DTLZ1_e20000} 
\caption{\label{fig:DTLZ1_macd_mana_20000} Pareto front after 20\,000 function evaluations, problem DTLZ1, 
$\mu$ARMOGA with population size 4 with crowding distance (left), and with the new proposed algorithm (right).}
\end{figure}

\begin{figure}[p]
\centering
\includegraphics[width=58mm]{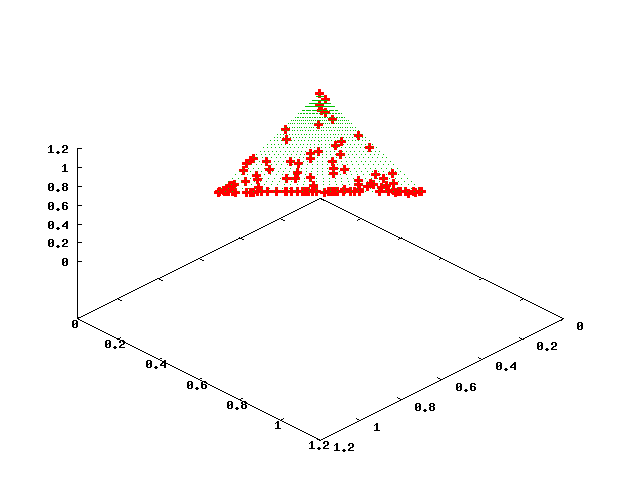}
\includegraphics[width=58mm]{emiarmoga4_DTLZ1_e200000} 
\caption{\label{fig:DTLZ1_macd_mana_200000} Pareto front after 200\,000 function evaluations, problem DTLZ1, 
$\mu$ARMOGA with population size 4 with crowding distance (left), and with the new proposed algorithm (right).}
\end{figure}

\begin{figure}[p]
\centering
\includegraphics[width=58mm]{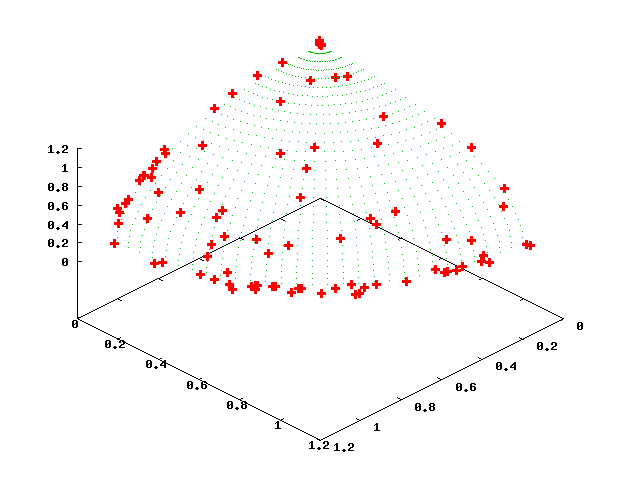}
\includegraphics[width=58mm]{emiarmoga4_DTLZ2_e4000} 
\caption{\label{fig:DTLZ2_macd_mana_4000} Pareto front after 4\,000 function evaluations, problem DTLZ2, 
$\mu$ARMOGA with population size 4 with crowding distance (left), and with the new proposed algorithm (right).}
\end{figure}

\begin{figure}[p]
\centering
\includegraphics[width=58mm]{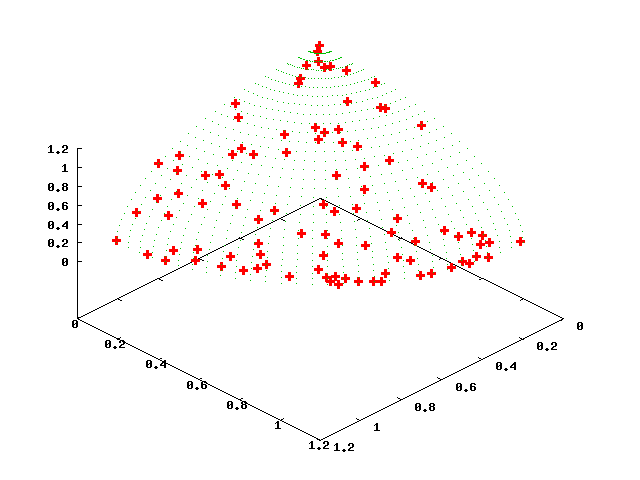}
\includegraphics[width=58mm]{emiarmoga4_DTLZ2_e20000} 
\caption{\label{fig:DTLZ2_macd_mana_20000} Pareto front after 20\,000 function evaluations, problem DTLZ2, 
$\mu$ARMOGA with population size 4 with crowding distance (left), and with the new proposed algorithm (right).}
\end{figure}

\begin{figure}[p]
\centering
\includegraphics[width=58mm]{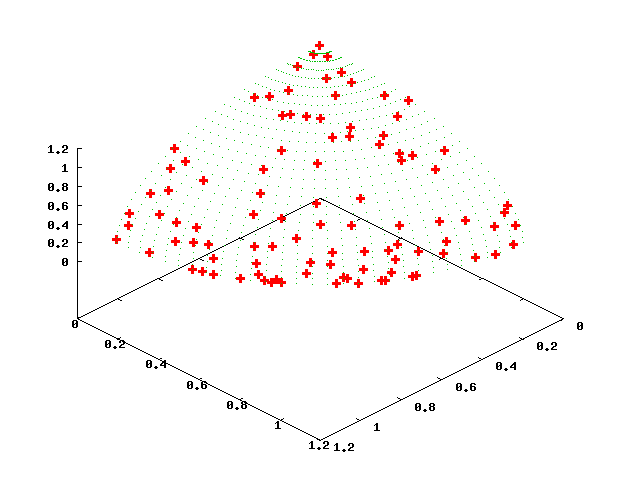}
\includegraphics[width=58mm]{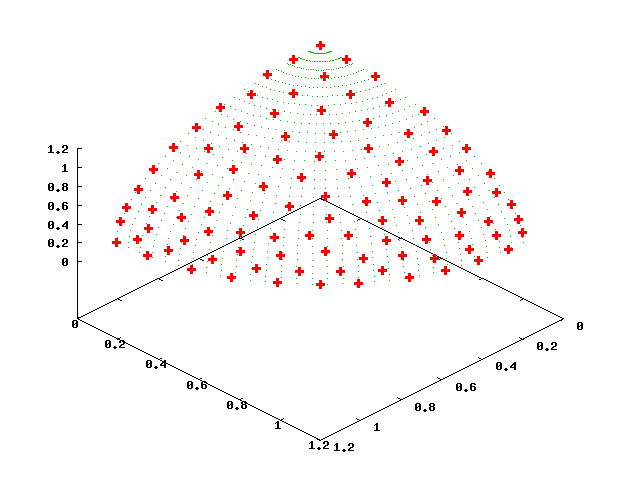} 
\caption{\label{fig:DTLZ2_macd_mana_200000} Pareto front after 200\,000 function evaluations, problem DTLZ2, 
$\mu$ARMOGA with population size 4 with crowding distance (left), and with the new proposed algorithm (right).}
\end{figure}

\begin{figure}[p]
\centering
\includegraphics[width=58mm]{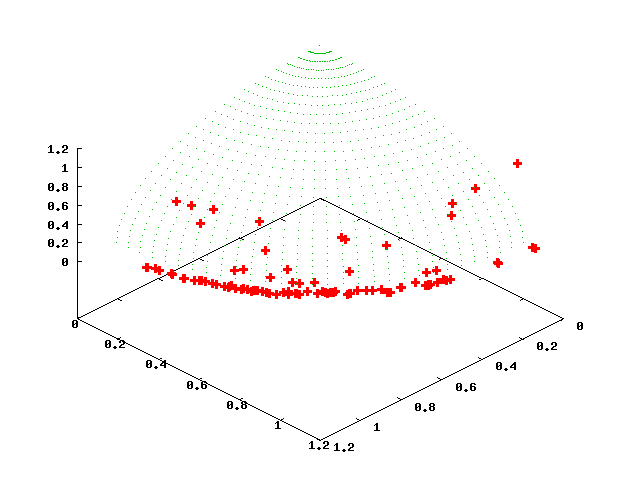}
\includegraphics[width=58mm]{emiarmoga4_DTLZ4_e4000} 
\caption{\label{fig:DTLZ4_macd_mana_4000} Pareto front after 4\,000 function evaluations, problem DTLZ4, 
$\mu$ARMOGA with population size 4 with crowding distance (left), and with the new proposed algorithm (right).}
\end{figure}

\begin{figure}[p]
\centering
\includegraphics[width=58mm]{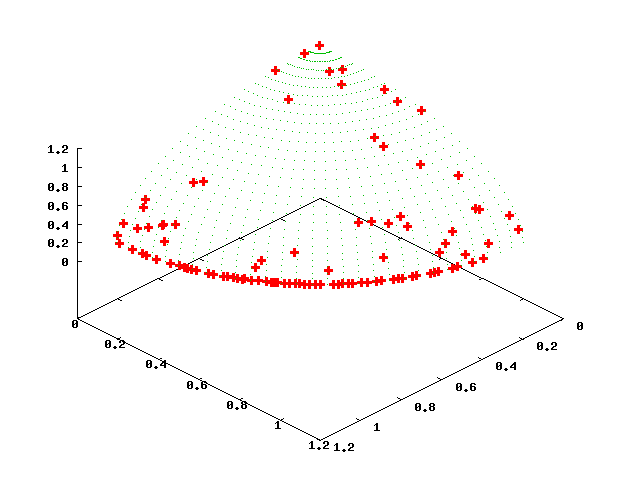}
\includegraphics[width=58mm]{emiarmoga4_DTLZ4_e20000} 
\caption{\label{fig:DTLZ4_macd_mana_20000} Pareto front after 20\,000 function evaluations, problem DTLZ4, 
$\mu$ARMOGA with population size 4 with crowding distance (left), and with the new proposed algorithm (right).}
\end{figure}

\begin{figure}[p]
\centering
\includegraphics[width=58mm]{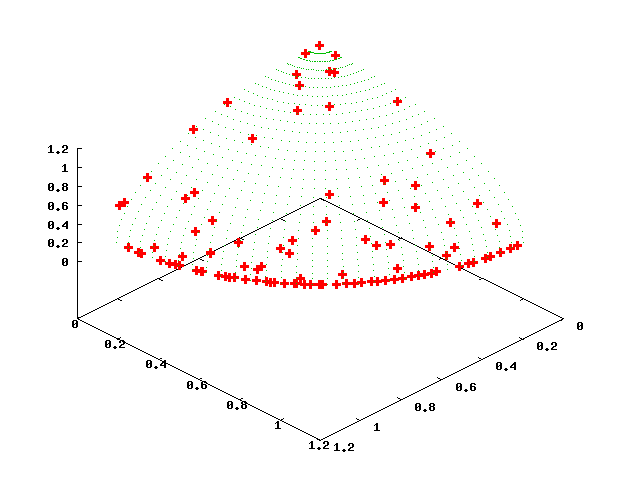}
\includegraphics[width=58mm]{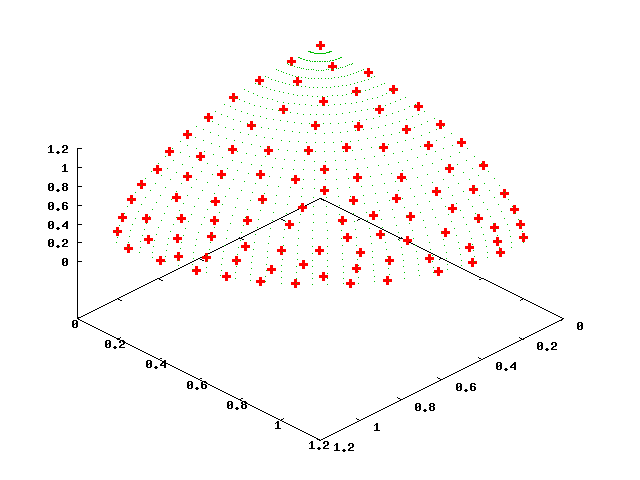} 
\caption{\label{fig:DTLZ4_macd_mana_200000} Pareto front after 200\,000 function evaluations, problem DTLZ4, 
$\mu$ARMOGA with population size 4 with crowding distance (left), and with the new proposed algorithm (right).}
\end{figure}

\begin{figure}[p]
\centering
\includegraphics[width=58mm]{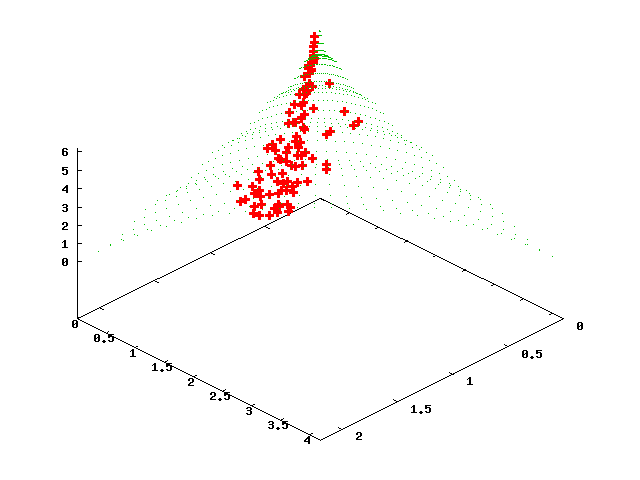}
\includegraphics[width=58mm]{emiarmoga4_WFG1_e200000}
\caption{\label{fig:WFG1_macd_mana_200000} Pareto front after 200\,000 function evaluations, problem WFG1, 
$\mu$ARMOGA with population size 4 with crowding distance (left), and with the new proposed algorithm (right).}
\end{figure}

\begin{figure}[p]
\centering
\includegraphics[width=58mm]{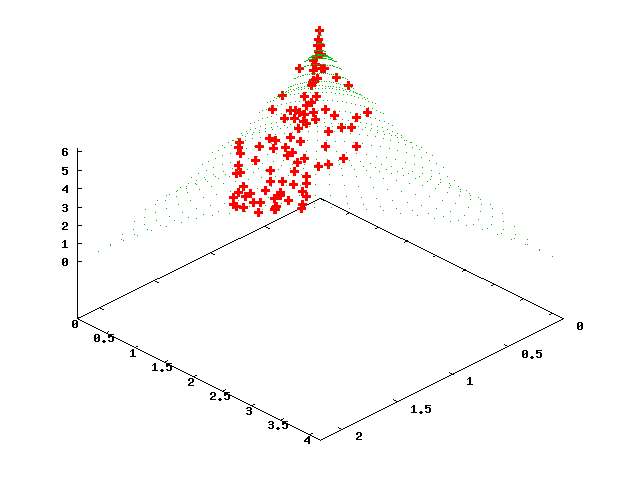}
\includegraphics[width=58mm]{emiarmoga4_WFG1_e1000000}
\caption{\label{fig:WFG1_macd_mana_1000000} Pareto front after 1\,000\,000 function evaluations, problem WFG1, 
$\mu$ARMOGA with population size 4 with crowding distance (left), and with the new proposed algorithm (right).}
\end{figure}

\begin{figure}[p]
\centering
\includegraphics[width=58mm]{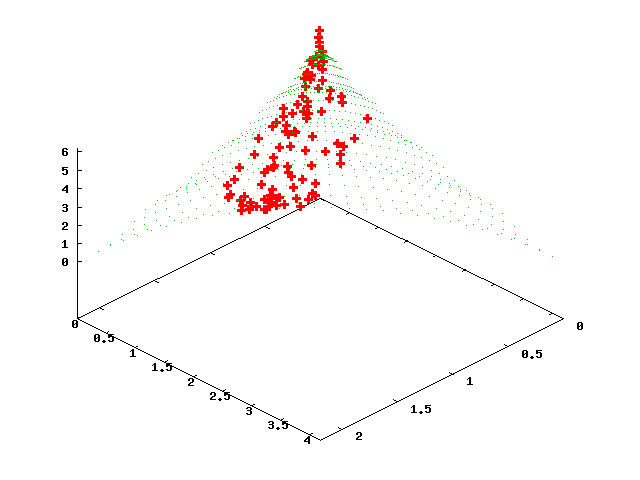}
\includegraphics[width=58mm]{emiarmoga4_WFG1_e2000000}
\caption{\label{fig:WFG1_macd_mana_2000000} Pareto front after 2\,000\,000 function evaluations, problem WFG1, 
$\mu$ARMOGA with population size 4 with crowding distance (left), and with the new proposed algorithm (right).}
\end{figure}

A set of experiments was run to compare $\mu$ARMOGA with crowding distance and $\mu$ARMOGA with the new archiving mechanism.
The population of four individuals was selected for the comparison.
Results for problems DTLZ1, DTLZ2, DTLZ4, and WFG1 are summarized in Tables \ref{tab:DTLZ1_macd}--\ref{tab:WFG1_mana} and 
Figures \ref{fig:DTLZ1_macd_mana_4000}--\ref{fig:WFG1_macd_mana_2000000}.

According to these experiments, the new archiving approach outperforms crowding distance in diversity as is clear 
from Figs.~\ref{fig:DTLZ1_macd_mana_4000}--\ref{fig:WFG1_macd_mana_2000000} and spacing metric in Tabs.~\ref{tab:DTLZ1_macd}--\ref{tab:WFG1_mana}.
While it also has very positive effect on the convergence of $\mu$ARMOGA to the exact Pareto front 
for problems DTLZ1, DTLZ2, and DTLZ4 (Tabs.~\ref{tab:DTLZ1_macd}--\ref{tab:DTLZ4_mana}),
both algorithms exhibit similar convergence for problem WFG1 (Tabs.~\ref{tab:WFG1_macd} and \ref{tab:WFG1_mana}).

\subsection{Summary}

As can be seen from above, $\mu$ARMOGA outperforms the other methods in distribution of individuals along the Pareto front,
and in many cases achieves the best convergence as well.
However, it should be noted 
that the default settings of the algorithms from PISA package is used, 
which may not be optimal for the test problems considered.

While IBEA offers exceptional convergence in some cases,
the distribution of individuals along the exact Pareto front is usually rather poor, 
with many individuals attached to limits of the objective space.
Our study confirms 
that the mechanism of crowding distance does not lead to uniform distribution of individuals along the Pareto front for more than two objectives.
The same result might be observed from the comparison of $\mu$ARMOGA using the two archiving mechanisms -- crowding distance and the new proposed technique
(Tables \ref{tab:DTLZ1_macd}--\ref{tab:WFG1_mana}, and Figures \ref{fig:DTLZ1_macd_mana_4000}--\ref{fig:WFG1_macd_mana_2000000}).
On the other hand, 
SPEA2 produces very uniform distribution of individuals comparable with $\mu$ARMOGA in some instances.

Concerning the number of evaluations,
DTLZ2 is the only problem,
for which only 4\,000 evaluations are sufficient to achieve reasonable convergence and distribution of individuals on the Pareto front by all algorithms.
On the other hand, for DTLZ1, even 40\,000 evaluations do not suffice to reach the true Pareto front for all seeds by any approach,
and results for 100\,000 and 200\,000 evaluations are added for a~reasonable comparison.
Even this large number of evaluations was not sufficient to reach the proximity of exact front in the case of problem WFG1,
and results for 1\,000\,000 and 2\,000\,000 evaluations are added.
For test functions DTLZ4 and WFG1, the new archiving mechanism is able to drive the evolution to regions, 
where the coverage of the Pareto front by individuals is sparse, 
and recover nice distribution of individuals along the Pareto set even for poorly chosen initial population.

To investigate the optimal distribution of the number of function evaluations between population size and number of generations for micro-evolution,
$\mu$ARMOGA is run with 4, 10 and 20 individuals for DTLZ1, DTLZ2, and DTLZ4.
According to our experiments, the performance of the algorithm is similar for all configurations with respect to spacing and convergence history and no strong dependence is revealed.
However, for problem DTLZ4, the method tends to produce more degenerated fronts with larger population (see Table~\ref{tab:degenerated}).
Additionally, population of 4 individuals leads to the best convergence metrics for problem DTLZ1,
and population of 20 individuals to the best converged front for problem DTLZ4.
Thus, using small populations and larger number of generations seems as the preferable approach.

\section{Conclusion}
\label{sec:conclusion}

The goal of our study is twofold:
(a) to develop a~new approach for selecting individuals to the Pareto archive;
(b) to explore the potential of using small population in evolutionary algorithms.

The main contribution of the paper is the presentation of a~new archiving mechanism.
Although its basic idea is rather simple and straightforward, the technique produces very promising results on all tested problems.
We are aware of the fact that the theoretical time complexity of the mechanism might be rather large (quadratic in the worst case).
However, our tests justify its usage, 
since the experimentally found complexity is much more favourable (approximately linear).
Moreover, it is intended to be used in combination with small population, 
for which such more elaborate selection mechanism is usually affordable.

The proposed selection mechanism was combined with $\mu$ARMOGA 
and is compared to other three state-of-the-art algorithms (NSGA-II, SPEA2, and IBEA) on four test problems.
We can conclude that $\mu$ARMOGA presents Pareto sets with the same or better distribution as SPEA2,
but usually with much better convergence to the exact front that is comparable with IBEA,
thus the best combining requirements on both convergence and distribution of individuals.
A considerable improvement is attained, using the new mechanism, in comparison with the version of $\mu$ARMOGA that uses crowding distance.
Clearly, $\mu$ARMOGA equipped with the new diversity mechanism is very promising and may be competitive with respect to other recent approaches.

Our experiments further support using small populations (up to 10 individuals), 
since runs with four individuals usually produces the best results.
It is well known that such small population can lead to rapid convergence.
However, in combination with the proposed archiving mechanism,
it also seems to be more robust with respect to an~initial population.

Regarding the history of convergence to the Pareto front,
in some cases as few as 4\,000 evaluations of objective function could be sufficient for some problems (DTLZ2),
while for other problems (multi-modal problem DTLZ1 or difficult WFG1), even 40\,000 evaluations may not be sufficient to approximate the true Pareto front,
and as many as 1\,000\,000 evaluations are needed for reasonable outcome.

\section*{Acknowledgement}
This research has been supported by the Development of Applied External Aerodynamics Program 
(Ministry of Education, Youth and Sports of the Czech Republic Grant MSM0001066901),
by research project AV0Z10190503 (Academy of Sciences of the Czech Republic),
and by grant IAA100760702 (Grant Agency of the Academy of Sciences of the Czech Republic).


\begin{thebibliography}{10}

\bibitem{Arakawa-1998-DAR}
{\sc Arakawa, M., and Hagiwara, I.}
\newblock Development of adaptive real range ({ARR}ange) genetic algorithms.
\newblock {\em JSME International Journal Series C, Mechanical systems, machine
  elements and manufacturing 41}, 4 (1998), 969--977.

\bibitem{Bleuler-2003-PIS}
{\sc Bleuler, S., Laumanns, M., Thiele, L., and Zitzler, E.}
\newblock {PISA} -- {A}~platform and programming language independent interface
  for search algorithms.
\newblock In {\em Evolutionary Multi-Criterion Optimization ({EMO 2003})\/}
  (Berlin, 2003), C.~M. Fonseca, P.~J. Fleming, E.~Zitzler, K.~Deb, and
  L.~Thiele, Eds., Lecture Notes in Computer Science, Springer, pp.~494--508.

\bibitem{Coello-2001-MGA}
{\sc Coello, C. A.~C., and Pulido, G.~T.}
\newblock A micro-genetic algorithm for multiobjective optimization.
\newblock In {\em Evolutionary Multi-Criterion Optimization ({EMO 2001})\/}
  (Berlin, 2001), E.~Zitzler, K.~Deb, L.~Thiele, C.~Coello~Coello, and
  D.~Corne, Eds., Lecture Notes in Computer Science, Springer, pp.~126--140.

\bibitem{Coello-2004-HMO}
{\sc Coello, C. A.~C., Pulido, G.~T., and Lechuga, M.~S.}
\newblock Handling multiple objectives with particle swarm optimization.
\newblock {\em IEEE Transactions on Evolutionary Computation 8}, 3 (2004),
  256--279.

\bibitem{Deb-2001-MOO}
{\sc Deb, K.}
\newblock {\em Multi-Objective Optimization using Evolutionary Algorithms}.
\newblock Wiley-Interscience Series in Systems and Optimization. John Wiley \&
  Sons, Chichester, 2001.

\bibitem{Deb-2000-FEN}
{\sc Deb, K., Agrawal, S., Pratap, A., and Meyarivan, T.}
\newblock A fast elitist non-dominated sorting genetic algorithm for
  multi-objective optimization: {NSGA-II}.
\newblock In {\em Parallel Problem Solving from Nature -- {PPSN VI}\/} (Berlin,
  2000), M.~Schoenauer, K.~Deb, G.~Rudolph, X.~Yao, E.~Lutton, J.~J. Merelo,
  and H.-P. Schwefel, Eds., Springer, pp.~849--858.

\bibitem{Deb-2003-FME}
{\sc Deb, K., Mohan, M., and Mishra, S.}
\newblock A fast multi-objective evolutionary algorithm for finding well-spread
  pareto-optimal solutions.
\newblock Tech. Rep. 2003002, KanGAL, Indian Institute of Technology Kanpur,
  India, 2003.

\bibitem{Deb-2002-SMO}
{\sc Deb, K., Thiele, L., Laumanns, M., and Zitzler, E.}
\newblock Scalable multi-objective optimization test problems.
\newblock In {\em Congress on Evolutionary Computation (CEC '02)\/}
  (Piscataway, New Jersey, 2002), IEEE Press, pp.~825--830.

\bibitem{Goldberg-1989-GAS}
{\sc Goldberg, D.~E.}
\newblock {\em Genetic Algorithms in Search, Optimization, and Machine
  Learning}.
\newblock {Addison-Wesley Professional}, January 1989.

\bibitem{Huband-2005-SMO}
{\sc Huband, S., Barone, L., While, R.~L., and Hingston, P.}
\newblock A scalable multi-objective test problem toolkit.
\newblock In {\em EMO\/} (2005), C.~A.~C. Coello, A.~H. Aguirre, and
  E.~Zitzler, Eds., vol.~3410 of {\em Lecture Notes in Computer Science},
  Springer, pp.~280--295.

\bibitem{Knowles-1999-PAE}
{\sc Knowles, J.~D., and Corne, D.~W.}
\newblock The {P}areto archived evolution strategy : A new baseline algorithm
  for {P}areto multiobjective optimisation.
\newblock In {\em Proceedings of the 1999 Congress on Evolutionary Computation
  (CEC'99)\/} (1999), vol.~1, pp.~98--105.

\bibitem{Knowles-2000-ANF}
{\sc Knowles, J.~D., and Corne, D.~W.}
\newblock Approximating the nondominated front using the {P}areto archived
  evolution strategy.
\newblock {\em Evolutionary Computation 8}, 2 (2000), 149--172.

\bibitem{Krishnakumar-1989-MGA}
{\sc Krishnakumar, K.}
\newblock Micro-genetic algorithms for stationary and non-stationary function
  optimization.
\newblock In {\em SPIE's Intelligent Control and Adaptive Systems Conference\/}
  (1989), G.~Rodriguez, Ed., Society of Photo-Optical Instrumentation Engineers
  (SPIE), pp.~289--296.

\bibitem{Oyama-2000-WDU}
{\sc Oyama, A.}
\newblock {\em Wing Design Using Evolutionary Algorithm}.
\newblock {PhD} thesis, Tohoku University, Japan, 2000.

\bibitem{Szollos-2008-AOM}
{\sc Sz\"{o}ll\"{o}s, A., \v{S}m\'{\i}d, M., and H\'{a}jek, J.}
\newblock Aerodynamic optimization via multi-objective micro-genetic algorithm
  with range adaptation, knowledge-based reinitialization, crowding and
  $\epsilon$-dominance.
\newblock {\em Advances in Engineering Software 40}, 6 (2009), 419--430.

\bibitem{Veldhuizen-1998-MEA}
{\sc Van~Veldhuizen, D.~A., and Lamont, G.~B.}
\newblock Multiobjective evolutionary algorithm research: A history and
  analysis.
\newblock Tech. Rep. TR-98-03, Air Force Institute of Technology,
  Wright-Patterson Air Force Base, Ohio, 2001.

\bibitem{Zitzler-2004-IBS}
{\sc Zitzler, E., and K{\"u}nzli, S.}
\newblock Indicator-based selection in multiobjective search.
\newblock In {\em Parallel Problem Solving from Nature (PPSN VIII)\/} (Berlin,
  2004), X.~Yao et~al., Eds., Lecture Notes in Computer Science, Springer,
  pp.~832--842.

\bibitem{Zitzler-2002-EMD}
{\sc Zitzler, E., Laumanns, M., and Thiele, L.}
\newblock {SPEA2}: Improving the strength pareto evolutionary algorithm for
  multiobjective optimization.
\newblock In {\em Evolutionary Methods for Design, Optimisation and Control
  with Application to Industrial Problems. Proceedings of the EUROGEN2001
  Conference\/} (Barcelona, 2002), K.~Giannakoglou, D.~Tsahalis, J.~Periaux,
  K.~Papaliliou, and T.~Fogarty, Eds., International Center for Numerical
  Methos in Engineering (CIMNE), pp.~95--100.

\bibitem{Zitzler-1999-MEA}
{\sc Zitzler, E., and Thiele, L.}
\newblock Multiobjective evolutionary algorithms: {A}~comparative case study
  and the strength pareto approach.
\newblock {\em IEEE Transactions on Evolutionary Computation 3}, 4 (1999),
  257--271.

\end{thebibliography}

\clearpage

\begin{table}[p]
\centering
\begin{tabular}{|c|c|c|c|}
	\hline
	metric  & NSGA-II & SPEA2 & IBEA \\
	\hline
	GD      & 3.61e+01 & 2.99e+01 & {\bf 1.80e+00}  \\
	TOL5    & 7.53e+01 & 5.92e+01 & {\bf 2.10e+00}  \\
	spacing & 2.47e+00 & 2.82e+00 & 2.14e+00  \\
	\hline
	\hline
	metric  & $\mu$ARMOGA(4) & $\mu$ARMOGA(10) & $\mu$ARMOGA(20)\\
	\hline
	GD      & 4.17e+00 & 1.13e+01 & 1.08e+01\\
	TOL5    & 5.70e+00 & 1.50e+01 & 1.50e+01\\
	spacing & {\bf 7.39e-01} & 7.40e-01 & 8.45e-01\\
	\hline
\end{tabular}
\caption{\label{tab:DTLZ1_4000}Problem DTLZ1, 4\,000 function evaluations}
\end{table}

\clearpage

\begin{table}[p]
\centering
\begin{tabular}{|c|c|c|c|}
	\hline
	metric  & NSGA-II & SPEA2 & IBEA \\
	\hline
	GD      & 4.96e+00 & 3.32e+00 & 5.91e-01 \\
	TOL5    & 8.52e+00 & 5.80e+00 & 6.64e-01 \\
	spacing & 1.50e+00 & 2.32e+00 & 2.65e+00 \\
	\hline
	\hline
	metric  & $\mu$ARMOGA(4) & $\mu$ARMOGA(10) & $\mu$ARMOGA(20)\\
	\hline
	GD      & {\bf 2.35e-01} & 6.10e+00 & 4.21e+00 \\
	TOL5    & {\bf 3.05e-01} & 7.44e+00 & 5.33e+00 \\
	spacing & {\bf 2.63e-01} & 9.85e-01 & 8.24e-01 \\
	\hline
\end{tabular}
\caption{\label{tab:DTLZ1_20000}Problem DTLZ1, 20\,000 function evaluations}
\end{table}

\clearpage

\begin{table}[p]
\centering
\begin{tabular}{|c|c|c|c|}
	\hline
	metric  & NSGA-II & SPEA2 & IBEA \\
	\hline
	GD      & 3.82e+00 & 2.31e+00 & 4.93e-01 \\
	TOL5    & 6.36e+00 & 3.71e+00 & 4.46e-01 \\
	spacing & 1.56e+00 & 1.71e+00 & 4.19e+00 \\
	\hline
	\hline
	metric  & $\mu$ARMOGA(4) & $\mu$ARMOGA(10) & $\mu$ARMOGA(20)\\
	\hline
	GD      & {\bf 3.61e-02} & 4.77e+00 & 3.02e+00 \\
	TOL5    & {\bf 5.04e-02} & 6.00e+00 & 3.79e+00 \\
	spacing & {\bf 1.39e-01} & 4.15e-01 & 1.70e-01 \\
	\hline
\end{tabular}
\caption{\label{tab:DTLZ1_40000}Problem DTLZ1, 40\,000 function evaluations}
\end{table}

\clearpage

\begin{table}[p]
\centering
\begin{tabular}{|c|c|c|c|}
	\hline
	metric  & NSGA-II & SPEA2 & IBEA \\
	\hline
	GD      & 2.33e+00 & 1.03e+00 & - \\
	TOL5    & 4.01e+00 & 2.01e+00 & - \\
	spacing & 1.51e+00 & 1.42e+00 & - \\
	\hline
	\hline
	metric  & $\mu$ARMOGA(4) & $\mu$ARMOGA(10) & $\mu$ARMOGA(20)\\
	\hline
	GD      & {\bf 1.88e-03} & 4.60e+00 & 1.91e+00 \\
	TOL5    & {\bf 4.48e-04} & 5.82e+00 & 2.41e+00 \\
	spacing & {\bf 8.50e-02} & 9.10e-02 & 1.51e-01 \\
	\hline
\end{tabular}
\caption{\label{tab:DTLZ1_100000}Problem DTLZ1, 100\,000 function evaluations}
\end{table}

\clearpage

\begin{table}[p]
\centering
\begin{tabular}{|c|c|c|c|}
	\hline
	metric  & NSGA-II & SPEA2 & IBEA \\
	\hline
	GD      & 2.26e+00 & 3.76e-01 & - \\
	TOL5    & 2.35e+00 & 5.64e-01 & - \\
	spacing & 1.43e+00 & 1.17e+00 & - \\
	\hline
	\hline
	metric  & $\mu$ARMOGA(4) & $\mu$ARMOGA(10) & $\mu$ARMOGA(20)\\
	\hline
	GD      & {\bf 1.54e-03} & 4.59e+00 & 2.08e+00 \\
	TOL5    & {\bf 1.74e-04} & 5.82e+00 & 2.30e+00 \\
	spacing & {\bf 7.82e-02} & 8.33e-02 & 5.19e-01 \\
	\hline
\end{tabular}
\caption{\label{tab:DTLZ1_200000}Problem DTLZ1, 200\,000 function evaluations}
\end{table}

\clearpage

\begin{table}[p]
\centering
\begin{tabular}{|c|c|c|c|}
	\hline
	metric  & NSGA-II & SPEA2 & IBEA \\
	\hline
	GD      & 1.07e-01 & 9.08e-02 & {\bf 4.96e-03} \\
	TOL5    & 2.51e-01 & 2.05e-01 & {\bf 7.49e-03} \\
	spacing & 6.40e-01 & 1.83e-01 & 7.44e-01 \\
	\hline
	\hline
	metric  & $\mu$ARMOGA(4) & $\mu$ARMOGA(10) & $\mu$ARMOGA(20)\\
	\hline
	GD      & 1.41e-02 & 1.77e-02 & 1.89e-02 \\
	TOL5    & 2.82e-02 & 3.77e-02 & 3.62e-02 \\
	spacing & {\bf 1.30e-01} & 1.64e-01 & 1.80e-01 \\
	\hline
\end{tabular}
\caption{\label{tab:DTLZ2_4000}Problem DTLZ2, 4\,000 function evaluations}
\end{table}

\clearpage

\begin{table}[p]
\centering
\begin{tabular}{|c|c|c|c|}
	\hline
	metric  & NSGA-II & SPEA2 & IBEA \\
	\hline
	GD      & 6.29e-02 & 3.50e-02 & {\bf 4.15e-04} \\
	TOL5    & 1.49e-01 & 7.53e-02 & {\bf 4.97e-04} \\
	spacing & 6.45e-01 & 1.40e-01 & 6.80e-01 \\
	\hline
	\hline
	metric  & $\mu$ARMOGA(4) & $\mu$ARMOGA(10) & $\mu$ARMOGA(20)\\
	\hline
	GD      & 2.59e-03 & 1.99e-03 & 2.94e-03 \\
	TOL5    & 4.01e-03 & 3.72e-03 & 3.29e-03 \\
	spacing & {\bf 6.94e-02} & 8.09e-02 & 8.75e-02 \\
	\hline
\end{tabular}
\caption{\label{tab:DTLZ2_20000}Problem DTLZ2, 20\,000 function evaluations}
\end{table}

\clearpage

\begin{table}[p]
\centering
\begin{tabular}{|c|c|c|c|}
	\hline
	metric  & NSGA-II & SPEA2 & IBEA \\
	\hline
	GD      & 4.64e-02 & 1.63e-02 & {\bf 4.94e-04} \\
	TOL5    & 1.06e-01 & 3.50e-02 & {\bf 2.71e-04} \\
	spacing & 6.14e-01 & 1.33e-01 & 6.82e-01 \\
	\hline
	\hline
	metric  & $\mu$ARMOGA(4) & $\mu$ARMOGA(10) & $\mu$ARMOGA(20)\\
	\hline
	GD      & 1.04e-03 & 1.32e-03 & 1.02e-03 \\
	TOL5    & 9.23e-04 & 1.02e-03 & 1.06e-03 \\
	spacing & {\bf 6.03e-02} & 6.71e-02 & 7.03e-02 \\
	\hline
\end{tabular}
\caption{\label{tab:DTLZ2_40000}Problem DTLZ2, 40\,000 function evaluations}
\end{table}

\clearpage

\begin{table}[p]
\centering
\begin{tabular}{|c|c|c|c|}
	\hline
	metric  & NSGA-II & SPEA2 & IBEA \\
	\hline
	GD      & 8.91e-02 & 1.31e-01 & 4.78e-03  \\
	TOL5    & 2.01e-01 & 2.81e-01 & 7.32e-03  \\
	spacing & 6.13e-01 & {\bf 1.80e-01} & 7.30e-01  \\
	\hline
	\hline
	metric  & $\mu$ARMOGA(4) & $\mu$ARMOGA(10) & $\mu$ARMOGA(20)\\
	\hline
	GD      & {\bf 1.87e-03} & 1.27e-02 & 5.45e-03 \\
	TOL5    & {\bf 3.42e-03} & 2.51e-02 & 1.23e-02 \\
	spacing & 8.48e-01 & 1.88e+00 & 2.09e+00 \\
	\hline
\end{tabular}
\caption{\label{tab:DTLZ4_4000}Problem DTLZ4, 4\,000 function evaluations}
\end{table}

\clearpage

\begin{table}[p]
\centering
\begin{tabular}{|c|c|c|c|}
	\hline
	metric  & NSGA-II & SPEA2 & IBEA  \\
	\hline
	GD      & 2.61e-02 & 3.92e-02 & 5.60e-04  \\
	TOL5    & 6.25e-02 & 9.73e-02 & 2.87e-04  \\
	spacing & 6.18e-01 & 1.21e-01 & 6.83e-01  \\
	\hline
	\hline
	metric  & $\mu$ARMOGA(4) & $\mu$ARMOGA(10) & $\mu$ARMOGA(20)\\
	\hline
	GD      & 4.34e-04 & 2.34e-04 & {\bf 3.55e-05} \\
	TOL5    & 1.43e-04 & 1.21e-04 & {\bf 3.40e-05} \\
	spacing & 1.46e-01 & {\bf 1.12e-01} & 3.45e-01 \\
	\hline
\end{tabular}
\caption{\label{tab:DTLZ4_20000}Problem DTLZ4, 20\,000 function evaluations}
\end{table}

\clearpage

\begin{table}[p]
\centering
\begin{tabular}{|c|c|c|c|}
	\hline
	metric  & NSGA-II & SPEA2 & IBEA \\
	\hline
	GD      & 1.61e-02 & 1.91e-02 & 1.75e-04 \\
	TOL5    & 2.94e-02 & 3.60e-02 & 1.39e-04 \\
	spacing & 6.42e-01 & 1.33e-01 & 6.94e-01 \\
	\hline
	\hline
	metric  & $\mu$ARMOGA(4) & $\mu$ARMOGA(10) & $\mu$ARMOGA(20)\\
	\hline
	GD      & 2.19e-04 & 3.81e-04 & {\bf 1.44e-05} \\
	TOL5    & 3.48e-05 & 2.78e-05 & {\bf 1.15e-05} \\
	spacing & 9.36e-02 & {\bf 8.01e-02} & 9.55e-02 \\
	\hline
\end{tabular}
\caption{\label{tab:DTLZ4_40000}Problem DTLZ4, 40\,000 function evaluations}
\end{table}

\clearpage

\begin{table}[p]
\centering
\begin{tabular}{|c|c|c|c|c|}
	\hline
	metric  & NSGA-II & SPEA2 & IBEA & $\mu$ARMOGA(4) \\
	\hline
	GD      & 8.86e-01 & 8.82e-01 & {\bf 7.47e-01} & 7.95e-01 \\
	TOL5    & 1.05e+00 & 1.04e+00 & {\bf 8.40e-01} & 9.27e-01 \\
	spacing & 4.90e-01 & {\bf 1.60e-01} & 5.40e-01 & 5.19e-01 \\
	\hline
\end{tabular}
\caption{\label{tab:WFG1_4000}Problem WFG1, 4\,000 function evaluations}
\end{table}

\begin{table}[p]
\centering
\begin{tabular}{|c|c|c|c|c|}
	\hline
	metric  & NSGA-II & SPEA2 & IBEA & $\mu$ARMOGA(4) \\
	\hline
	GD      & 9.01e-01 & 8.95e-01 & 7.67e-01 & {\bf 4.97e-01} \\
	TOL5    & 1.06e+00 & 1.05e+00 & 8.62e-01 & {\bf 5.92e-01} \\
	spacing & 4.88e-01 & {\bf 1.67e-01} & 5.66e-01 & 2.82e-01 \\
	\hline
\end{tabular}
\caption{\label{tab:WFG1_20000}Problem WFG1, 20\,000 function evaluations}
\end{table}

\begin{table}[p]
\centering
\begin{tabular}{|c|c|c|c|c|}
	\hline
	metric  & NSGA-II & SPEA2 & IBEA & $\mu$ARMOGA(4) \\
	\hline
	GD      & 8.86e-01 & 8.82e-01 & 7.47e-01 & {\bf 3.77e-01} \\
	TOL5    & 1.05e+00 & 1.04e+00 & 8.40e-01 & {\bf 4.56e-01} \\
	spacing & 4.90e-01 & {\bf 1.60e-01} & 5.40e-01 & 2.24e-01 \\
	\hline
\end{tabular}
\caption{\label{tab:WFG1_40000}Problem WFG1, 40\,000 function evaluations}
\end{table}

\begin{table}[p]
\centering
\begin{tabular}{|c|c|c|c|c|}
	\hline
	metric  & NSGA-II & SPEA2 & IBEA & $\mu$ARMOGA(4) \\
	\hline
	GD      & 8.80e-01 & 8.75e-01 & 7.26e-01 & {\bf 2.51e-01} \\
	TOL5    & 1.05e+00 & 1.04e+00 & 8.15e-01 & {\bf 3.15e-01} \\
	spacing & 4.94e-01 & {\bf 1.62e-01} & 5.29e-01 & 1.96e-01 \\
	\hline
\end{tabular}
\caption{\label{tab:WFG1_100000}Problem WFG1, 100\,000 function evaluations}
\end{table}

\clearpage

\begin{table}[p]
\centering
\begin{tabular}{|c|c|c|c|c|}
	\hline
	metric  & NSGA-II & SPEA2 & IBEA & $\mu$ARMOGA(4) \\
	\hline
	GD      & 8.75e-01 & 8.79e-01 & 7.10e-01 & {\bf 1.75e-01} \\
	TOL5    & 1.04e+00 & 1.04e+00 & 7.99e-01 & {\bf 2.29e-01} \\
	spacing & 4.98e-01 & {\bf 1.66e-01} & 5.47e-01 & 1.86e-01 \\
	\hline
\end{tabular}
\caption{\label{tab:WFG1_200000}Problem WFG1, 200\,000 function evaluations}
\end{table}

\clearpage

\begin{table}[p]
\centering
\begin{tabular}{|c|c|c|c|c|}
	\hline
	metric  & NSGA-II & SPEA2 & IBEA & $\mu$ARMOGA(4) \\
	\hline
	GD      & 8.58e-01 & 8.71e-01 & 6.81e-01 & {\bf 6.53e-02} \\
	TOL5    & 1.03e+00 & 1.03e+00 & 7.66e-01 & {\bf 1.14e-01} \\
	spacing & 4.72e-01 & 1.74e-01 & 5.24e-01 & {\bf 1.53e-01} \\
	\hline
\end{tabular}
\caption{\label{tab:WFG1_1000000}Problem WFG1, 1\,000\,000 function evaluations}
\end{table}

\clearpage

\begin{table}[p]
\centering
\begin{tabular}{|c|c|c|c|c|}
	\hline
	metric  & NSGA-II & SPEA2 & IBEA & $\mu$ARMOGA(4) \\
	\hline
	GD      & 8.58e-01 & 8.74e-01 & 6.69e-01 & {\bf 4.79e-02} \\
	TOL5    & 1.03e+00 & 1.03e+00 & 7.52e-01 & {\bf 9.69e-02} \\
	spacing & 5.06e-01 & 1.89e-01 & 5.35e-01 & {\bf 1.54e-01} \\
	\hline
\end{tabular}
\caption{\label{tab:WFG1_2000000}Problem WFG1, 2\,000\,000 function evaluations}
\end{table}

\clearpage

\begin{table}[p]
\centering
\begin{tabular}{|c|c|c|c|}
	\hline
	metric  & $\mu$ARMOGA(4) & $\mu$ARMOGA(10) & $\mu$ARMOGA(20)\\
	\hline
	DTLZ1   & 0.98 & 1.32 & 1.32 \\
	DTLZ2   & 0.99 & 1.42 & 1.47 \\
	DTLZ4   & 1.05 & 1.35 & 1.35 \\
	WFG1    & 1.08 & n/a  & n/a  \\
	\hline
\end{tabular}
\caption{\label{tab:updates}Average number of updates necessary after the addition of an~individual into archive, archive size limit 100.}
\end{table}

\begin{table}[p]
\centering
\begin{tabular}{|c|c|c|c|c|c|c|}
	\hline
	archive size  & 20 & 50 & 100 & 200 & 500 & 1\,000 \\
	\hline
	DTLZ4   & 1.06  & 1.05 & 1.06 & 1.06 & 1.02 & 0.94 \\
	\hline
\end{tabular}
\caption{\label{tab:asize_updates}Average number of updates necessary after the addition of an~individual into archive for problem DTLZ4, variable archive size.}
\end{table}

\begin{table}[p]
\centering
\begin{tabular}{|c|c|c|c|}
	\hline
	metric  & NSGA-II & SPEA2 & IBEA \\
	\hline
	DTLZ1   & -/-/-/-/- & -/-/-/-/- & -/2/11/20/20 \\
	DTLZ2   & -/-/- & -/-/- & -/-/- \\
	DTLZ4   & 11/8/5 & 11/11/11 & 6/6/6 \\
	WFG1    & -/-/-/-/-/-/- & -/-/-/-/-/-/- & -/-/-/-/-/-/- \\
	\hline
	\hline
	metric  & $\mu$ARMOGA(4) & $\mu$ARMOGA(10) & $\mu$ARMOGA(20)\\
	\hline
	DTLZ1   & -/-/-/-/- & -/-/-/-/- & -/-/-/-/- \\
	DTLZ2   & -/-/- & -/-/- & -/-/- \\
	DTLZ4   & -/-/- & 4/4/4 & 4/3/3 \\
	WFG1    & -/-/-/-/-/-/- & n/a & n/a \\
	\hline
\end{tabular}
\caption{\label{tab:degenerated}Number of degenerated Pareto fronts for 20 seeds. 
Number at 4\,000/20\,000/40\,000(/100\,000/200\,000(/1\,000\,000/2\,000\,000)) function evaluations.}
\end{table}

\begin{table}[p]
\centering
\begin{tabular}{|c|c|c|c|c|c|}
	\hline
	evaluations  & 4\,000 & 20\,000 & 40\,000 & 100\,000 & 200\,000 \\
	\hline
	GD           & 2.13e+01 & 2.28e+00 & 1.15e+00 & 2.41e+00 & 3.72e-01 \\
	TOL5         & 4.02e+01 & 1.42e+00 & 2.07e-01 & 3.99e-02 & 1.22e-02 \\
	spacing      & 3.34e+00 & 2.34e+00 & 2.73e+00 & 3.29e+00 & 2.60e+00 \\
	\hline
        degenerated  & 0        & 0        & 0        & 0        & 0        \\
	\hline
\end{tabular}
\caption{\label{tab:DTLZ1_macd}$\mu$ARMOGA with crowding distance,
four individuals, problem DTLZ1, average for twenty seeds.}
\end{table}

\begin{table}[p]
\centering
\begin{tabular}{|c|c|c|c|c|c|}
	\hline
	evaluations  & 4\,000 & 20\,000 & 40\,000 & 100\,000 & 200\,000 \\
	\hline
	GD           & 4.17e+00 & 2.35e-01 & 3.61e-02 & 1.88e-03 & 1.54e-03 \\
	TOL5         & 5.70e+00 & 3.05e-01 & 5.04e-02 & 4.48e-04 & 1.74e-04 \\
	spacing      & 7.39e-01 & 2.63e-01 & 1.39e-01 & 8.50e-02 & 7.82e-02 \\
	\hline
        degenerated  & 0        & 0        & 0        & 0        & 0 \\
	\hline
\end{tabular}
\caption{\label{tab:DTLZ1_mana}$\mu$ARMOGA with the new proposed archiving algorithm,
four individuals, problem DTLZ1, average for twenty seeds.}
\end{table}

\begin{table}[p]
\centering
\begin{tabular}{|c|c|c|c|c|c|}
	\hline
	evaluations  & 4\,000 & 20\,000 & 40\,000 & 100\,000 & 200\,000 \\
	\hline
	GD                     & 5.36e-02 & 2.54e-03 & 1.77e-03 & 2.16e-03 & 1.99e-03 \\
	TOL5                   & 7.92e-02 & 2.92e-03 & 2.83e-03 & 3.58e-03 & 2.83e-01 \\
	spacing                & 5.74e-01 & 5.60e-01 & 5.47e-01 & 5.49e-01 & 5.41e-01 \\
	\hline
        degenerated            & 0            & 0            & 0            & 0            & 0 \\
	\hline
\end{tabular}
\caption{\label{tab:DTLZ2_macd}$\mu$ARMOGA with crowding distance,
four individuals, problem DTLZ2, average for twenty seeds.}
\end{table}

\begin{table}[p]
\centering
\begin{tabular}{|c|c|c|c|c|c|}
	\hline
	evaluations  & 4\,000 & 20\,000 & 40\,000 & 100\,000 & 200\,000 \\
	\hline
	GD                     & 1.41e-02 & 2.59e-03 & 1.04e-03 & 3.45e-04 & 1.81e-04 \\
	TOL5                   & 2.82e-02 & 4.01e-03 & 9.23e-04 & 1.09e-04 & 2.39e-05 \\
	spacing                & 1.30e-01 & 6.94e-02 & 6.03e-02 & 4.94e-02 & 4.46e-02 \\
	\hline
        degenerated            & 0            & 0            & 0            & 0            & 0 \\
	\hline
\end{tabular}
\caption{\label{tab:DTLZ2_mana}$\mu$ARMOGA with the new proposed archiving algorithm,
four individuals, problem DTLZ2, average for twenty seeds.}
\end{table}

\begin{table}[p]
\centering
\begin{tabular}{|c|c|c|c|c|c|}
	\hline
	evaluations  & 4\,000 & 20\,000 & 40\,000 & 100\,000 & 200\,000 \\
	\hline
	GD                     & 5.46e+01 & 2.35e+01 & 1.18e+01 & 1.68e+01 & 1.48e+01 \\
	TOL5                   & 1.46e+02 & 1.96e+01 & 1.92e-01 & 1.70e-01 & 1.67e-01 \\
	spacing                & 3.40e+00 & 6.50e+00 & 5.09e+00 & 4.73e+00 & 4.86e+00 \\
	\hline
        degenerated            & 2            & 3            & 3            & 3            & 3 \\
	\hline
\end{tabular}
\caption{\label{tab:DTLZ4_macd}$\mu$ARMOGA with crowding distance,
four individuals, problem DTLZ4, average for twenty seeds.}
\end{table}

\begin{table}[p]
\centering
\begin{tabular}{|c|c|c|c|c|c|}
	\hline
	evaluations  & 4\,000 & 20\,000 & 40\,000 & 100\,000 & 200\,000 \\
	\hline
	GD                     & 1.87e-03 & 4.34e-04 & 2.19e-04 & 9.32e-05 & 3.39e-05 \\
	TOL5                   & 3.42e-03 & 1.43e-04 & 3.48e-05 & 1.34e-05 & 8.82e-06 \\
	spacing                & 8.48e-01 & 1.46e-01 & 9.36e-02 & 7.00e-02 & 6.35e-02 \\
	\hline
        degenerated            & 0            & 0            & 0            & 0            & 0 \\
	\hline
\end{tabular}
\caption{\label{tab:DTLZ4_mana}$\mu$ARMOGA with the new proposed archiving algorithm,
four individuals, problem DTLZ4, average for twenty seeds.}
\end{table}

\begin{table}[p]
\centering
\begin{tabular}{|c|c|c|c|c|c|c|c|}
	\hline
	evaluations  & 4\,000 & 20\,000 & 40\,000 & 100\,000 & 200\,000 & 1\,000\,000 & 2\,000\,000 \\
	\hline
	GD           & 9.01e-01 & 3.67e-01 & 2.81e-01 & 1.93e-01 & 1.27e-01 & 4.09e-02 & 3.93e-02 \\
	TOL5         & 1.21e+00 & 4.26e-01 & 3.32e-01 & 2.37e-01 & 1.66e-01 & 8.01e-02 & 8.51e-02 \\
	spacing      & 1.46e+00 & 7.03e-01 & 6.85e-01 & 6.17e-01 & 6.24e-01 & 5.83e-01 & 5.28e-01 \\
	\hline
        degenerated            & 0            & 0            & 0            & 0            & 0            & 0            & 0 \\
	\hline
\end{tabular}
\caption{\label{tab:WFG1_macd}$\mu$ARMOGA with crowding distance,
four individuals, problem WFG1, average for twenty seeds.}
\end{table}

\begin{table}[p]
\centering
\begin{tabular}{|c|c|c|c|c|c|c|c|}
	\hline
	evaluations  & 4\,000 & 20\,000 & 40\,000 & 100\,000 & 200\,000 & 1\,000\,000 & 2\,000\,000 \\
	\hline
	GD           & 7.95e-01 & 4.97e-01 & 3.77e-01 & 2.51e-01 & 1.75e-01 & 6.53e-02 & 4.79e-02 \\
	TOL5         & 9.27e-01 & 5.92e-01 & 4.56e-01 & 3.15e-01 & 2.29e-01 & 1.14e-01 & 9.69e-02 \\
	spacing      & 5.19e-01 & 2.82e-01 & 2.24e-01 & 1.96e-01 & 1.86e-01 & 1.53e-01 & 1.54e-01 \\
	\hline
        degenerated            & 0            & 0            & 0            & 0            & 0            & 0            & 0 \\
	\hline
\end{tabular}
\caption{\label{tab:WFG1_mana}$\mu$ARMOGA with the new proposed archiving algorithm,
four individuals, problem WFG1, average for twenty seeds.}
\end{table}

\clearpage

\clearpage

\end{document}